\documentclass[11pt,english,a4paper]{amsart}

\usepackage{babel}
\usepackage[left=2.5cm,right=2.5cm,top=2.5cm,bottom=2.5cm]{geometry}
\usepackage[latin1]{inputenc}
\usepackage{amsmath}
\usepackage{amsfonts}
\usepackage{dsfont}
\usepackage{amssymb}
\usepackage{graphicx,latexsym,xspace,epsfig}
\usepackage[dvipsnames]{xcolor}
\usepackage{caption}
\usepackage{subcaption}
\usepackage{color}
\usepackage{ulem}
\usepackage{amsthm}
\usepackage{natbib} 
\usepackage{numprint} 
\usepackage{setspace} 
\usepackage[pdftex]{hyperref}
\usepackage[]{algorithm2e}
\usepackage{enumitem} 

\usepackage{placeins}
\usepackage{afterpage}

\numberwithin{equation}{section}
\newtheorem{theoreme}{Theorem}[section]
\newtheorem{prop}[theoreme]{Proposition}
\newtheorem{cor}[theoreme]{Corollary}
\newtheorem{lemme}[theoreme]{Lemma}
\newtheorem{definition}{Definition}[section]

\newtheorem{assumption}{Assumption}
 
\numberwithin{equation}{section}

\renewcommand{\vec}[1]{{\boldsymbol{#1}}}

\newcommand{\mat}[1]{{\boldsymbol{#1}}}

\title[Conditional extreme risk measures from heavy-tailed elliptical random vectors]{Estimation of conditional extreme risk measures from heavy-tailed elliptical random vectors}

\author{A. Usseglio-Carleve}
\address{Universit\'e de Lyon, Universit\'e Lyon 1,  Institut Camille Jordan ICJ UMR 5208 CNRS}
\email{usseglio@math.univ-lyon1.fr }   

\begin{document}

\begin{abstract}
In this work, we focus on some conditional extreme risk measures estimation for elliptical random vectors. In a previous paper, we proposed a methodology to approximate extreme quantiles, based on two extremal parameters. We thus propose some estimators for these parameters, and study their consistency and asymptotic normality in the case of heavy-tailed distributions. Thereafter, from these parameters, we construct extreme conditional quantiles estimators, and give some conditions that ensure consistency and asymptotic normality. Using recent results on the asymptotic relationship between quantiles and other risk measures, we deduce estimators for extreme conditional $L_p-$quantiles and Haezendonck-Goovaerts risk measures. Under similar conditions, consistency and asymptotic normality are provided. In order to test the effectiveness of our estimators, we propose a simulation study. A financial data example is also proposed.
\end{abstract}

\maketitle

\small{
\noindent{Keywords:}~\textit{\keywords{Elliptical distribution; Extreme quantiles; Extreme value theory; Haezendonck-Goovaerts risk measures; Heavy-tailed distributions; $L_p-$quantiles.}}}

\section{Introduction} \label{Introduction}

In many fields such as finance or actuarial science, quantile, or Value-at-Risk (see~\cite{linsmeier}) is a recognized tool for risk measurement. In~\cite{koenker}, quantile is seen as minimum of an asymmetric loss function. However, Value-at-Risk, or VaR, has some disadvantages, such as that of not being a coherent measure in the sense of~\cite{artzner}. These limits have led many authors to use alternative risk measures. \\ On the basis of Koenker's approach,~\cite{newey} proposed another measure called expectile, which has since been widely studied (see for example~\cite{sobotka} or more recently~\cite{daouia}) and applied (\cite{taylor} and~\cite{cai}). Later,~\cite{mquantiles} introduced M-quantiles, a family of measures minimizing an asymmetric loss function, and~\cite{chen} focused on asymmetric power functions to define $L_p-$quantiles. The cases $p=1$ and $p=2$ correspond respectively to the quantile and expectile. Recently,~\cite{bignozzi} provided some results concerning $L_p-$quantiles for Student distributions, and have shown that closed formula are difficult to obtain in the general case. \\
In parallel,~\cite{artzner} introduced the Tail-Value-at-Risk as an alternative to Value-at-Risk, and this risk measure subsequently had many applications (see e.g.~\cite{marceau}). Moreover, TVaR belongs to a larger family of risk measures called Haezendonck-Goovaerts risk measures and introduced in~\cite{haezendonck},~\cite{goovaerts} and~\cite{tangyang}. In the same way as $L_p-$quantiles, we do not have an explicit formula in the general case. \\
However, for a heavy-tailed random variable,~\cite{girardlp} proved that $L_p-$quantile and $L_1-$quantile (or quantile) are asymptotically proportional. Then, as proposed in~\cite{daouia}, an estimator of a $L_p-$quantile may be deduced from a suitable estimator of the quantile, for extreme levels. In the same spirit,~\cite{tangyang} provided a similar asymptotic relationship between a subclass of Haezendonck-Goovaerts risk measures and quantiles. Finally, all these risk measures we introduced may be estimated through a quantile estimation in an asymptotic setting. \\
Extreme quantiles estimation is a very active area of research. In recent years, we can give many examples~: \cite{gardesweibull} focused on Weibull tail distributions, \cite{elmethni} proposed a study for heavy and light tailed distributions, \cite{gong} was interested in functions of dependent variables, and \cite{devalk} provided a methodology for high quantiles estimation. The question of extreme conditional quantiles estimation has also been explored in \cite{wang} in a regression framework. However, \cite{article} and \cite{article2} have shown that the regression setting may lead to a poor estimation of extreme measures in the case of elliptical distributions.  Elliptical distributions, introduced in~\cite{kelker}, aim to generalize the gaussian distribution, i.e to define symmetric distributions with different properties, such as a heavy tail. This is why elliptical distributions are more and more used in finance (see for example~\cite{owen} or~\cite{xiao}). \\
For all these reasons, we consider, in this paper, an elliptical random vector $\vec{Z}=(\vec{X},Y)$ with the consistency property (in the sense of~\cite{kano}), where $\vec{X} \in \mathbb{R}^N$, $Y \in \mathbb{R}$, and propose to estimate some extreme quantiles (and deduce $L_p-$quantiles and Haezendonck-Goovaerts risk measures) of $Y|\vec{X}=\vec{x}$, i.e. of a component conditionnally to the others. In order to improve the conditional quantile estimation, we proposed in~\cite{article} a methodology based on two extremal parameters, and the unconditional quantile of $Y$. Indeed, if we denote $F_{Y|\vec{x}}^{-1}(\alpha)$ the quantile of level $\alpha$ of $Y|\vec{X}=\vec{x}$, the latter is asymptotically equivalent to a quantile of $Y$ ($F_{Y}^{-1}$ will be the quantile function of $Y$), in the following manner~:

\begin{equation}
F_{Y|\vec{x}}^{-1}(\alpha) \underset{\alpha \rightarrow 1}{\sim} F_{Y}^{-1} \left( \delta(\alpha, \eta, \ell) \right),
\label{relationquantiles}
\end{equation}
where $\delta$ is a known function (detailed later) depending on $\alpha$ and two parameters $\eta$ and $\ell$ called extremal parameters. One can notice that Equation~\eqref{relationquantiles} may only holds under the consistency property of $\vec{Z}$. \cite{article} has also shown that extremal parameters do not exist for some consistent elliptical distributions (see e.g. the Laplace distribution). \\ In this paper, the goal will be in a first time to give a sufficient condition on $\vec{Z}$ that ensures the existence of $\eta$ and $\ell$. This is why a regularly varying assumption is done. After having proved their existence, estimators for the parameters $\eta$ and $\ell$ are proposed, and therefore for extreme conditional quantiles. \\
The paper is organized as follows. Section~\ref{prelim} provides some definitions and properties of elliptical distributions, including the extremal parameters introduced in~\cite{article}. A particular interest is given to consistent elliptical distributions. Section~\ref{parametres} is devoted to extremal parameters $\eta$ and $\ell$. Under a regularly varying assumption, their existence is proved, and estimators are proposed. By adding some conditions, consistency and asymptotic normality results are given.  In Section~\ref{quantiles}, we use the results of Section~\ref{parametres} to introduce some estimators of extreme quantiles, and give consistency and asymptotic normality results. The asymptotic relationships between $L_p-$quantiles and quantiles recalled in Section~\ref{lpquantiles} allow us to give extreme $L_p-$quantiles estimators. The same approach is proposed for extreme Haezendonck-Goovaerts risk measures. In order to analyze the efficiency of our estimators, we propose a simulation study in Section~\ref{simstudy}, and a real data example in Section~\ref{real}.

\section{Preliminaries} \label{prelim}

In this section, we first recall some classical results on elliptical distributions. We consider a $d-$dimensional vector $\vec{Z}$ from an elliptical distribution with parameters $\vec{\mu} \in \mathbb{R}^d$ and $\mat{\Sigma} \in \mathbb{R}^{d \times d}$. Then the density of $\vec{Z}$, if it exists, is given by~:
\begin{equation}
\frac{c_d}{|\mat{\Sigma}|^{\frac{1}{2}}} g_d \left( \left(\vec{z}-\vec{\mu} \right)^T \mat{\Sigma}^{-1}\left(\vec{z}-\vec{\mu} \right) \right).
\label{densite}
\end{equation}
$c_d$ and $g_d$ will respectively be called normalization coefficient and generator of $\vec{Z}$. \cite{cambanis} gives another way to characterize an elliptical distribution, through the following stochastic representation~:
\begin{equation}
\vec{Z} \overset{d}{=} \vec{\mu}+R \mat{\Lambda} \vec{U}^{(d)},
\label{repcambanis}
\end{equation}
where $\mat{\Lambda} \mat{\Lambda}^T=\mat{\Sigma}$, $\vec{U}^{(d)}$ is a $d-$dimensional random vector uniformly distributed on the unit sphere of dimension $d$, and $R$ is a non-negative random variable independent of $\vec{U}^{(d)}$. $R$ is called radius of $\vec{Z}$. In the following, the radius must have a particular shape. Indeed, \cite{cambanis1979} and \cite{kano} propose a representation for some particular elliptical distributions. Let us consider $\left( \vec{Z}_d \right)_{d \in \mathbb{N}^*}$ a family of elliptical distributions of dimension $d$. Then $\left( \vec{Z}_d \right)_{d \in \mathbb{N}^*}$ possesses the consistency property if it admits the following representation for all $d \in \mathbb{N}^*$~:
\begin{equation}
\vec{Z}_d \overset{d}{=} \vec{\mu}+\chi_d \xi \mat{\Lambda} \vec{U}^{(d)},
\label{eqkano}
\end{equation}
where $\chi_d$ is the square root of a $\chi^2$ distribution with $d$ degrees of freedom, $\xi$ is a non-negative random variable which does not depend on $d$, and $\chi_d$, $\xi$ and $\vec{U}^{(d)}$ are mutually independent. In \cite{kano}, such elliptical distributions are said consistent, have the advantage of being stable by linear combinations (combining Theorem 2.16 of~\cite{fang1990} and Theorem 1 in~\cite{kano}), and allow us to define elliptical random fields (see, e.g.,~\cite{opitz}). In the following, we focus on consistent elliptical distributions, and take the notation
\begin{equation}
R_d=\chi_d \xi.
\label{loiR}
\end{equation}
For the sake of clarity, we will say that a random variable with stochastic representation~\eqref{eqkano} is $(\xi,d)-$elliptical with parameters $\vec{\mu}$ and $\mat{\Sigma}$. Using this terminology, the purpose of the paper is as follows. Let $\vec{Z}=(\vec{X},Y) \in \mathbb{R}^{N+1}$ be a $(\xi,N+1)-$elliptical random vector with parameters $\vec{\mu}$ and $\mat{\Sigma}$, where $\vec{X} \in \mathbb{R}^N$ and $Y \in \mathbb{R}$. Consistency property of $\vec{Z}$ implies that $\vec{X}$ and $Y$ are respectively $(\xi,N)-$ and $(\xi,1)-$elliptical distributions with parameters $\vec{\mu_X} \in \mathbb{R}^N$, $\mat{\Sigma_X} \in \mathbb{R}^{N \times N}$ and $\mu_Y \in \mathbb{R}$, $\Sigma_Y \in \mathbb{R}$. We also denote $\vec{\Sigma}_{\vec{X}Y}$ the covariance vector between $\vec{X}$ and $Y$. The aim is thus to provide a predictor for the quantile of the conditional distribution $Y|\vec{X}=\vec{x}$. According to Theorem 7 of~\cite{frahm}, such a distribution is still elliptical, with a radius $R^*$ different from $R$ in the general case. In particular, we have~:
\begin{equation}
\{Y|\vec{X}=\vec{x} \} \overset{d}{=} \mu_{Y|\vec{X}}+\sigma_{Y|\vec{X}} R^* U^{(1)},
\end{equation}
where $\mu_{Y|\vec{X}}=\mu_Y+\vec{\Sigma}_{\vec{X}Y}^T \mat{\Sigma_{X}}^{-1}(\vec{x}-\vec{\mu_X})$ and $\sigma_{Y|\vec{X}}^2=\Sigma_{Y}-\vec{\Sigma}_{\vec{X}Y}^T \mat{\Sigma_{X}}^{-1} \vec{\Sigma}_{\vec{X}Y}$.
Then, denoting $\Phi_{R^*}(t)=\mathbb{P} \left( R^* U^{(1)} \leq t \right)$, and using the translation equivariance and positive homogeneity of elliptical quantiles (see~\cite{mcneil2015quantitative}), conditional quantiles of $Y|\vec{X}=\vec{x}$ may be expressed as~:
\begin{equation}
q_{\alpha}(Y|\vec{X}=\vec{x})= \mu_{Y|\vec{X}}+\sigma_{Y|\vec{X}} \Phi_{R^*}^{-1}(\alpha),
\label{eqcondquant}
\end{equation}
where $\alpha \in ]0,1[$. Thus, in order to give a good prediction of $q_{\alpha}(Y|\vec{X}=\vec{x})$, we need to estimate the conditional function $\Phi_{R^*}^{-1}$. Unfortunately, when we have a data set $\vec{X}_1,...,\vec{X}_n$, we only observe the unconditional distribution of $\vec{X}$. This is why, in~\cite{article}, we have given a predictor for conditional quantiles, based solely on the unconditional c.d.f $\Phi_R(t)=\mathbb{P} \left( R_1 U^{(1)} \leq t \right)$. This approximation is based on two parameters $\eta \in \mathbb{R}$ and $0< \ell< +\infty$ such that~:
\begin{equation}
\underset{t \rightarrow \infty}{\lim} \text{ } \frac{\bar{\Phi}_{R^*}(t)}{\bar{\Phi}_{R}(t^{\eta})}=\ell.
\label{limite}
\end{equation}
Table~\ref{tablecoeffs} gives some examples of coefficients $\eta$ and $\ell$ for classical elliptical distributions. 
\begin{table}[!h]
\begin{center}
\begin{tabular}{lccc}
  Distribution & $\eta$ & $\ell$ \\
  \hline
  Gaussian & $1$  & $1$ \\ [2mm]
  Student, $\nu >0 $ & $\frac{N}{\nu}+1$ & $\frac{\Gamma \left(\frac{\nu+N+1}{2} \right) \Gamma \left(\frac{\nu}{2} \right)}{\Gamma \left(\frac{\nu+N}{2}\right) \Gamma \left(\frac{\nu+1}{2} \right)} \left(1+\frac{M(\vec{x})}{\nu} \right)^{\frac{N+\nu}{2}} \frac{\nu^{\frac{N}{2}+1}}{\nu+N}$ \\ [2mm]
  UGM & $1$ & $\frac{\rm \min(\theta_1, \ldots,\theta_n)^N \exp \left\{ -\frac{\rm \min(\theta_1, \ldots,\theta_n)^2}{2} M(\vec{x}) \right\}  }{\sum \limits_{k=1}^n \pi_k \theta_k^N \exp \left( -\frac{\theta_k^2}{2} M(\vec{x}) \right)}$ \\ [2mm]
  Slash, $a >0 $ & $\frac{N}{a}+1$ & $\frac{\Gamma \left(\frac{N+1+a}{2} \right) M(\vec{x})^{\frac{N+a}{2}} }{\Gamma \left( \frac{N+a}{2} \right) (N+a) \chi_{N+a}^2 \left(M(\vec{x}) \right) 2^{\frac{a}{2}-1} \Gamma \left( \frac{1+a}{2} \right) }$ \\
  \hline
\end{tabular}
\end{center}
\caption{Coefficients $\eta$ and $\ell$ for classical distributions, where $M(\vec{x})=(\vec{x}-\vec{\mu_X})^{\top} \mat{\Sigma_{X}}^{-1}(\vec{x}-\vec{\mu_X})$.}
\label{tablecoeffs}
\end{table}
However, we have shown in~\cite{article} that such parameters not always exist for all elliptical distribution (see, e.g, Laplace distribution). In a first time, we can wonder in which setting these parameters exist. We thus consider the following assumption, that will ensures the existence of $\eta$ and $\ell$.

\begin{assumption}[Second order regular variations]
We assume that there exist a function $A$ such that $A(t) \rightarrow 0$ as $t \rightarrow +\infty$, and
\begin{equation}
\underset{t \rightarrow +\infty}{\lim} \frac{\frac{\Phi_{R}^{-1} \left(1-\frac{1}{\omega t} \right)}{\Phi_{R}^{-1} \left(1-\frac{1}{t} \right)}-\omega^{\gamma}}{A(t)}=\omega^{\gamma} \frac{\omega^{\rho}-1}{\rho},
\label{eqhyp2}
\end{equation}
where $\gamma>0$ and $\rho < 0$.
\label{hyp2}
\end{assumption}

This assumption is widespread in literature of extreme quantiles (see, e.g, \cite{daouia}). A first consequence is that $\Phi_R$, or equivalently $F_{R_1}$ is attracted to the maximum domain of Pareto-type distributions with tail index $\gamma$. Furthermore, it entails $\Phi_{R}^{-1}(1-1/t) \sim c_1 t^{\gamma}$, or equivalently $\bar{\Phi}_{R}(t) \sim c_2 t^{-\frac{1}{\gamma}}$ as $t \rightarrow +\infty$ (see~\cite{dehaan2006}). As example, Student distribution satisfies Assumption~\ref{hyp2}. The following lemma provides some results concerning asymptotic equivalences.

\begin{lemme}[Regular variation properties]
Under Assumption~\ref{hyp2}, we get the following regular variations properties~:
\begin{enumerate}[label=(\roman*)]
\item The random variable $\xi$ satisfies
\begin{equation}
\bar{F}_{\xi}(t) \underset{t \rightarrow +\infty}{\sim} \lambda t^{-\frac{1}{\gamma}}, \lambda \in \mathbb{R}.
\label{lemme1_0}
\end{equation}
\item For all $d \in \mathbb{N}^*$, the random variable $R_d=\chi_d \xi$ is attracted to the maximum domain of Pareto-type distribution with tail index $\gamma$, and
\begin{equation}
\bar{F}_{R_d}(t) \underset{t \rightarrow +\infty}{\sim} 2^{\frac{1}{2 \gamma}} \frac{\Gamma \left( \frac{d+\gamma^{-1}}{2} \right)}{\Gamma \left( \frac{d}{2} \right)} \bar{F}_{\xi}(t) \underset{t \rightarrow +\infty}{\sim} 2^{\frac{1}{2 \gamma}} \frac{\Gamma \left( \frac{d+\gamma^{-1}}{2} \right)}{\Gamma \left( \frac{d}{2} \right)} \lambda t^{-\frac{1}{\gamma}}, \lambda \in \mathbb{R}.
\label{lemme1_1}
\end{equation}
\item For all $\eta>0$, $d \in \mathbb{N}^*$,  
\begin{equation}
\frac{f_{R_d}(t)}{f_{R_1}(t^{\eta})}\underset{t \rightarrow +\infty}{\sim} \frac{\sqrt{\pi} \Gamma \left( \frac{d+\gamma^{-1}}{2} \right)}{\Gamma \left( \frac{d}{2} \right) \Gamma \left( \frac{1+\gamma^{-1}}{2} \right)} t^{(\eta-1)(\gamma^{-1}+1)}.
\label{lemme1_3}
\end{equation}
\end{enumerate}
\label{lemmeprops}
\end{lemme}
These results will be usefull throughout the paper, and especially in the following result which proves the existence of our parameters.

\begin{prop}[Existence of extremal parameters]
Under Assumption~\ref{hyp2}, parameters $\eta$ and $\ell$ exist, and are expressed~:
\begin{equation}
\left \{
\begin{array}{cc}
       \eta=&  1 +\gamma N\\
       \ell=& \frac{\Gamma \left( \frac{N+\gamma^{-1}+1}{2} \right)}{\Gamma \left( \frac{\gamma^{-1}+1}{2} \right)} \frac{\gamma^{-1} \pi^{-\frac{N}{2}}}{(N+\gamma^{-1}) c_N g_N\left(M(\vec{x}) \right)} \\
\end{array}
\right..
\label{eqetal}
\end{equation}
\label{propetal}
\end{prop}
One can notice that $\eta$ is only related to the tail index $\gamma$, and not to the covariate vector $\vec{x}$, while $\ell$ is depending on $c_N g_N \left( M(\vec{x}) \right)$. In the next, we thus denote rather $\ell(\vec{x})$, in order to emphasize the role played by the covariate vector $\vec{x}$. We can now give the following predictor for $q_{\alpha}\left( Y|\vec{X}=\vec{x} \right)$~: 
\begin{equation}
q_{\alpha \uparrow}(Y|\vec{X}=\vec{x})= \mu_{Y|\vec{X}}+\sigma_{Y|\vec{X}} \left[ \Phi_{R}^{-1}\left\{ 1-\frac{1}{\frac{\ell(\vec{x})}{1-\alpha}+2(1-\ell(\vec{x}))} \right\} \right] ^{{1}/{\eta}}.
\label{ecqe}
\end{equation}
From there, we have proved in Theorem 7 of~\cite{article} that $q_{\alpha \uparrow}(Y|\vec{X}=\vec{x})$ and $q_{\alpha}(Y|\vec{X}=\vec{x})$ were asymptoticaly equivalent as $\alpha \rightarrow 1$, i.e 
\begin{equation}
q_{\alpha \uparrow}(Y|\vec{X}=\vec{x}) \underset{\alpha \rightarrow 1}{\sim} q_{\alpha}(Y|\vec{X}=\vec{x}).
\label{equivalencephirstar}
\end{equation}
A similar equivalence has been easily deduced for $\alpha \rightarrow 0$, using the symmetry properties of elliptical distributions. In this paper, we focus on the case $\alpha \rightarrow 1$, case $\alpha \rightarrow 0$ being easily deduced.  
In Section~\ref{parametres}, we propose some estimators for extremal parameters $\eta$ and $\ell(\vec{x})$. Before that, we need to do a little simplification. Indeed, Equation~\eqref{ecqe} shows that the extreme quantile estimation requires  the prior estimation of quantities $\mu_{Y|\vec{X}}$ and $\sigma_{Y|\vec{X}}$. These quantities may be easily estimated by the method of moments or fixed-point algorithm (c.f p.66 of~\cite{frahm}). In a spatial setting, even if the variable $Y$ is not observed, a stationarity assumption on the random field makes it possible to estimate these values (see~\cite{cressie}). Furthermore, the speed of convergence of these methods is higher than those of the estimators we propose in this paper, and therefore do not interfere in the asymptotic results. This is why, in the following, we suppose that $\mu_{Y|\vec{X}}$, $\sigma_{Y|\vec{X}}$, and therefore $\vec{\mu_X}$, $\mat{\Sigma_X}$ are known. Then, it remains to estimate $\eta$, $\ell(\vec{x})$ and $\Phi_{R^*}^{-1}$. Section~\ref{parametres} focuses on $\eta$ and $\ell(\vec{x})$, while Section~\ref{quantiles} deals with $\Phi_{R^*}^{-1}$.

\section{Extremal coefficients estimation} \label{parametres}

In this section, the aim is to estimate the extremal parameters $\eta$ and $\ell(\vec{x})$ conditionaly to the covariates vector $\vec{X}=\vec{x}$. For that purpose, we consider a random sample $\vec{X}_1,...,\vec{X}_n$ independent and identically distributed from an $(\xi,N)-$elliptical vector with the same distribution as $\vec{X}$, and denote $M(\vec{x})=\left( \vec{x}-\vec{\mu_X} \right)^T \mat{\Sigma_X}^{-1} \left( \vec{x}-\vec{\mu_X} \right)$. The aim is then to give two suitable estimator $\hat{\eta}$ and $\hat{\ell}(\vec{x})$, respectively for $\eta$ and $\ell(\vec{x})$.

\subsection{Estimation of $\eta$}

We notice that coefficient $\eta$ is directly related on the tail index $\gamma$. Then, using a suitable estimator of $\gamma$, we easily deduce $\eta$. There are several estimators widespread in the literature. As examples, \cite{pickands}, \cite{schultze} or \cite{kratz} provide some estimators for $\gamma$. In the following, we use the Hill estimator, introduced in \cite{hill}~:
\begin{equation}
\hat{\gamma}_{k_n}=\frac{1}{k_n} \sum \limits_{i=1}^{k_n} \ln \left( \frac{W_{[i]}}{W_{[k_n+1]}} \right),
\label{hillest}
\end{equation}
where $W_{[1]} \geq \hdots \geq W_{[k_n+1]} \geq \hdots \geq W_{[n]}$ and $k_n=o(n)$ such that $k_n \rightarrow +\infty$ as $n \rightarrow +\infty$. In this context, the statistic $W$ may be~:
\begin{itemize}
\item The first (or indifferently any) component of the reduced centered covariate vector $\mat{\Lambda_X}^{-1} \left(\vec{X}-\vec{\mu_X} \right)$, where $\mat{\Lambda_X}^T \mat{\Lambda_X}=\mat{\Sigma_X}$. This approach works well, but we do not use all available data.
\item The Mahalanobis norm $\sqrt{(\vec{X}-\vec{\mu_X})^T \mat{\Sigma_{X}}^{-1}(\vec{X}-\vec{\mu_X})}$. This approach has the advantage of using all available data.
\end{itemize}
Indeed, according to Theorem 2 of~\cite{hashorvanorm}, the two last quantities both admit $\gamma$ as tail index. \\
In the following we will use the one-component approach, since the asymptotic results we give are valid under Assumption~\ref{hyp2}, applied to the univariate c.d.f $\Phi_R$. Moreover, numerical comparisons seem show that the second approach does not significantly improve the estimation of the parameters. Main properties of $\hat{\gamma}_{k_n}$ may be found in \cite{dehaanhill}. Under second order condition given in Assumption~\ref{hyp2}, \cite{dehaan2006} proved the following asymptotic normality for $\hat{\gamma}_{k_n}$.
\begin{equation}
\sqrt{k_n} \left( \hat{\gamma}_{k_n}-\gamma \right) \underset{n \rightarrow +\infty}{\rightarrow} \mathcal{N} \left( \frac{\lambda}{1-\rho}, \gamma^2 \right),
\label{eqnormgammahat}
\end{equation}
where $\lambda=\underset{n \rightarrow +\infty}{\lim} \sqrt{k_n} A \left(  \frac{n}{k_n} \right)$ and $k_n=o(n)$ such that $k_n \rightarrow +\infty$ as $n \rightarrow +\infty$. Then, using Proposition~\ref{propetal} and Equation~\eqref{hillest}, we define the following estimator for $\eta$.

\begin{definition}[Estimator of $\eta$]
We define $\hat{\eta}_{k_n}$ as
\begin{equation}
\hat{\eta}_{k_n}=\frac{N}{k_n} \sum \limits_{i=1}^{k_n} \ln \left( \frac{W_{[i]}}{W_{[k_n+1]}} \right)+1.
\label{eqetahat}
\end{equation}
\label{defetahat}
\end{definition}
As an affine transformation of Hill estimator, asymptotic normality of $\hat{\eta}_{k_n}$ is obvious. In order to simplify the next results, we suppose $\lambda=0$ in what follows.

\begin{prop}[Asymptotic normality of $\hat{\eta}_{k_n}$]
Under Assumption~\ref{hyp2}, and if $\underset{n \rightarrow +\infty}{\lim} \sqrt{k_n} A \left(  \frac{n}{k_n} \right)=0$, then
\begin{equation}
\sqrt{k_n} \left( \hat{\eta}_{k_n}-\eta \right) \underset{n \rightarrow +\infty}{\rightarrow} \mathcal{N} \left( 0, N^2 \gamma^2 \right).
\end{equation}
\label{propetahat}
\end{prop}

\subsection{Estimation of $\ell(\vec{x})$} The form of $\ell(\vec{x})$, given in Proposition~\ref{propetal}, leads to a more complicated estimation. Indeed, $\ell(\vec{x})$ is related on both $\gamma$ and $c_N g_N\left(M(\vec{x}) \right)$. Our estimator for $\gamma$ is given in Equation~\eqref{hillest}. Concerning $c_N g_N\left(M(\vec{x}) \right)$, we propose a kernel estimator. Class of kernel estimators, introduced in \cite{parzen}, makes it possible to estimate probability densities. Then, the following lemma will be usefull for the construction of our estimator. This result comes from p.108 of \cite{johnson}.
\begin{lemme}
The Mahalanobis distance $M(\vec{X})=\left( \vec{X}-\vec{\mu_X} \right)^T \mat{\Sigma_{X}}^{-1} \left( \vec{X}-\vec{\mu_X} \right)$ has density~:
\begin{equation}
f_{M(\vec{X})}(t)=\frac{\pi^{\frac{N}{2}}}{\Gamma \left( \frac{N}{2} \right)} x^{\frac{N}{2}-1} c_N g_N(t).
\label{densitemaha}
\end{equation}
\label{lemmemaha}
\end{lemme}  
Using Lemma~\ref{lemmemaha}, we introduced a kernel estimator $\hat{g}_{h_n}$ for $c_N g_N\left(M(\vec{x}) \right)$.
\begin{definition}[Generator estimator]
We define $\hat{g}_{h_n}$ as
\begin{equation}
\hat{g}_{h_n}=\frac{M(\vec{x})^{1-\frac{N}{2}} \Gamma \left( \frac{N}{2} \right) }{\pi^{\frac{N}{2}}} \hat{f}_{M(\vec{X})}(M(\vec{x}))=\frac{M(\vec{x})^{1-\frac{N}{2}} \Gamma \left( \frac{N}{2} \right) }{\pi^{\frac{N}{2}} nh_n} \sum \limits_{i=1}^n K \left( \frac{M(\vec{x})-\left(\vec{X}_i-\vec{\mu_X} \right)^T \mat{\Sigma_{X}}^{-1}\left(\vec{X}_i-\vec{\mu_X} \right)}{h_n} \right),
\label{eqghat}
\end{equation}
where the kernel $K$ fills some conditions given in \cite{parzen} and bandwith $h_n$ verifies $h_n \rightarrow 0$ and $n h_n \rightarrow +\infty$ as $n \rightarrow +\infty$.
\label{defghat}
\end{definition}
\cite{parzen} provided the asymptotic normality for kernel estimators. We first define some assumptions concerning $K$ and $g_N$ needed for the next results.
\begin{itemize}
\item $(K1)$ : $K$ is compactly supported on $[-1,1]$ and bounded. In addition, $\int_{\mathbb{R}} K(u) du=1$, $K(u)=K(-u) \text{ } \forall u \in \mathbb{R}$ and $\int_{\mathbb{R}} u^2 K(u) du \neq 0$.
\item $(K2)$ : In the neighborhood of $M(\vec{x})$, $g_N$ is bounded and twice continuously differentiable with bounded derivatives.
\end{itemize} 
The following results may be found in~\cite{li2007}. Under conditions $(K1)-(K2)$, it may be proved that~:
\begin{equation}
\left \{
\begin{array}{cc}
       \mathbb{E}\left[  \hat{f}_{M(\vec{X})}(M(\vec{x})) \right]- f_{M(\vec{X})}(M(\vec{x})) &=  O(h_n^2)\\
       \mathbb{V}ar \left[ \hat{f}_{M(\vec{X})}(M(\vec{x})) \right] &= O\left( \frac{1}{n h_n} \right) \\
\end{array}
\right..
\label{grandtaus}
\end{equation}
By adding the condition $nh_n^5 \rightarrow 0$ as $n \rightarrow +\infty$, we also obtain the asymptotic normality~:
\begin{equation}
\sqrt{nh_n} \left( \hat{f}_{M(\vec{X})}(M(\vec{x})) - f_{M(\vec{X})}(M(\vec{x})) \right) \underset{n \rightarrow +\infty}{\rightarrow} \mathcal{N} \left( 0, f_{M(\vec{X})}(M(\vec{x})) \int K(u)^2 du \right).
\label{eqnormfhat}
\end{equation}
Using the previous results given above, the following asymptotic normality for $\hat{g}_{h_n}$ is easily deduced.
\begin{prop}[Asymptotic normality of generator estimator]
Under conditions $(K1)-(K2)$, and taking a sequence $h_n$ such that $h_n \rightarrow 0$, $n h_n \rightarrow +\infty$ and $n h_n^5 \rightarrow 0$ as $n \rightarrow +\infty$, then the following relationship holds~:
\begin{equation}
\sqrt{nh_n} \left( \hat{g}_{h_n}-c_N g_N\left(M(\vec{x}) \right) \right) \underset{n \rightarrow +\infty}{\rightarrow} \mathcal{N} \left( 0, \frac{M(\vec{x})^{1-\frac{N}{2}} \Gamma \left( \frac{N}{2} \right)}{\pi^{\frac{N}{2}}} c_N g_N\left(M(\vec{x}) \right) \int K(u)^2 du \right).
\label{eqnormghat}
\end{equation}
\label{propnormghat}
\end{prop}
Replacing $\gamma$ by $\hat{\gamma}$ and $c_N g_N\left(M(\vec{x}) \right)$ by $\hat{g}_{h_n}$ in Equation~\eqref{eqetal}, we are now able to provide an estimator $\hat{\ell}(\vec{x})$ for $\ell(\vec{x})$, in the following definition. Furthermore, under Assumption~\ref{hyp2}, we give the asymptotic normality of $\hat{\ell}(\vec{x})$.

\begin{definition}[Estimator of $\ell(\vec{x})$]
We define $\hat{\ell}_{k_n,h_n}(\vec{x})$ as~:
\begin{equation}
\hat{\ell}_{k_n,h_n}(\vec{x})=\frac{\Gamma \left( \frac{N+\hat{\gamma}_{k_n}^{-1}+1}{2} \right)}{\Gamma \left( \frac{\hat{\gamma}_{k_n}^{-1}+1}{2} \right)} \frac{\hat{\gamma}_{k_n}^{-1} \pi^{-\frac{N}{2}}}{\left(N+\hat{\gamma}_{k_n}^{-1} \right) \hat{g}_{h_n}}.
\label{eqlhat}
\end{equation}
where $\hat{\gamma}_{k_n}$ and $\hat{g}_{h_n}$ are respectively given in Equations~\eqref{hillest} and \eqref{eqghat}.
\label{proplhat}
\end{definition}

\begin{prop}
Under Assumption~\ref{hyp2}, conditions $(K1)-(K2)$ and if $\underset{n \rightarrow +\infty}{\lim} \sqrt{k_n} A \left( \frac{n}{k_n} \right)=0$, the following asymptotic relationships hold~:
\begin{enumerate}[label=(\roman*)]
\item If $nh_n/k_n \underset{n \rightarrow +\infty}{\rightarrow} +\infty$ and $\sqrt{k_n} h_n^2 \underset{n \rightarrow +\infty}{\rightarrow} 0$, then
\begin{equation}
\sqrt{k_n} \left( \hat{\ell}_{k_n,h_n}(\vec{x})-\ell(\vec{x}) \right) \underset{n \rightarrow +\infty}{\rightarrow} \mathcal{N} \left( 0, V_1 \left(\gamma,c_N g_N\left(M(\vec{x}) \right)\right) \right).
\label{eqnorml}
\end{equation}
\item If $nh_n/k_n \underset{n \rightarrow +\infty}{\rightarrow} 0$ and $nh_n^5 \underset{n \rightarrow +\infty}{\rightarrow} 0$, then
\begin{equation}
\sqrt{nh_n} \left( \hat{\ell}_{k_n,h_n}(\vec{x})-\ell(\vec{x}) \right) \underset{n \rightarrow +\infty}{\rightarrow} \mathcal{N} \left( 0, V_2 \left(\gamma,c_N g_N\left(M(\vec{x}) \right) \right) \right),
\label{eqnorml2}
\end{equation}
\end{enumerate}
where ($\Psi$ is the digamma function (see p.258 of~\cite{handbook}))
\begin{equation}
\left \{
\begin{array}{cc}
      V_1(\gamma,c_N g_N\left(M(\vec{x}) \right)) =& \frac{\pi^{-N} \gamma^2}{c_N^2g_N\left(M(\vec{x}) \right)^2}  \frac{\Gamma \left( \frac{N+\gamma^{-1}+1}{2} \right)^2}{\Gamma \left( \frac{\gamma^{-1}+1}{2} \right)^2} \left[ \frac{\Psi\left( \frac{\gamma^{-1}+1}{2} \right)-\Psi\left( \frac{N+\gamma^{-1}+1}{2} \right)}{2 \gamma^2 \left( N \gamma+1 \right)}-\frac{N}{\left( N \gamma+1 \right)^2} \right]^2 \\
      V_2(\gamma,c_N g_N\left(M(\vec{x}) \right)) =& \frac{ \Gamma \left( \frac{N}{2} \right)}{M(\vec{x})^{\frac{N}{2}-1} \pi^{\frac{N}{2}}} c_N g_N \left(M(\vec{x}) \right) \int K(u)^2 du \left[ \frac{\Gamma \left( \frac{N+\gamma^{-1}+1}{2} \right)}{\Gamma \left( \frac{\gamma^{-1}+1}{2} \right)} \frac{\gamma^{-1} \pi^{-\frac{N}{2}}}{\left( N+\gamma^{-1} \right) c_N^2 g_N\left(M(\vec{x}) \right)^2} \right]^2 \\
\end{array}
\right..
\label{variances}
\end{equation}
\label{propnorml}
\end{prop}

We have the asymptotic normality for our estimators $\hat{\eta_{k_n}}$ and $\hat{\ell}_{k_n,h_n}(\vec{x})$. The next proposition gives the joint distribution according to the asymptotic relations between $k_n$ and $h_n$. The proof derives from delta method. 

\begin{prop}
Under Assumption~\ref{hyp2}, conditions $(K1)-(K2)$ and if $\sqrt{k_n} A \left( \frac{n}{k_n} \right) \rightarrow 0$ as $n \rightarrow +\infty$, then the following asymptotic relationships hold~:
\begin{enumerate}[label=(\roman*)]
\item If $nh_n/k_n \underset{n \rightarrow +\infty}{\rightarrow} 0$ and $nh_n^5 \underset{n \rightarrow +\infty}{\rightarrow} 0$, then
\begin{equation}
\sqrt{n h_n} \begin{pmatrix}
   \hat{\ell}_{k_n,h_n}(\vec{x})-\ell(\vec{x})  \\
   \hat{\eta}_{k_n}-\eta  
\end{pmatrix} \underset{n \rightarrow +\infty}{\rightarrow} \mathcal{N} \left( \begin{pmatrix}
   0  \\
   0  
\end{pmatrix}, \begin{pmatrix}
   V_2(\gamma,c_N g_N\left(M(\vec{x}) \right)) & 0 \\
   0 & 0 
\end{pmatrix} \right),
\label{eqloijointe1}
\end{equation}
where $V_2(\gamma,c_N g_N\left(M(\vec{x}) \right))$ is given in Equation~\eqref{variances}.
\item If $nh_n/k_n \underset{n \rightarrow +\infty}{\rightarrow} +\infty$ and $\sqrt{k_n} h_n^2 \underset{n \rightarrow +\infty}{\rightarrow} 0$, then
\begin{equation}
\sqrt{k_n} \begin{pmatrix}
   \hat{\ell}_{k_n,h_n}(\vec{x})-\ell(\vec{x})  \\
   \hat{\eta}_{k_n}-\eta  
\end{pmatrix} \underset{n \rightarrow +\infty}{\rightarrow} \mathcal{N} \left( \begin{pmatrix}
   0  \\
   0  
\end{pmatrix}, \begin{pmatrix}
   V_1(\gamma,c_N g_N\left(M(\vec{x}) \right)) & -N \gamma \sqrt{V_1(\gamma,c_N g_N\left(M(\vec{x}) \right))} \\
   -N \gamma \sqrt{V_1(\gamma,c_N g_N\left(M(\vec{x}) \right))} & N^2 \gamma^2 
\end{pmatrix} \right),
\label{eqloijointe2}
\end{equation}
where $V_1(\gamma,c_N g_N\left(M(\vec{x}) \right))$ is given in Equation~\eqref{variances}.
\end{enumerate}
\label{proploijointe}
\end{prop}
Using the previous results, we propose, in Section~\ref{quantiles}, some estimators of extreme conditional quantiles based on $\hat{\ell}_{k_n,h_n}(\vec{x})$ and $\hat{\eta}_{k_n}$.

\section{Extreme quantiles estimation} \label{quantiles}

In this section, we propose some estimators of extreme quantiles $q_{\alpha_n} \left( Y|\vec{X}=\vec{x} \right)$, for a sequence $\alpha_n \rightarrow 1$ as $n \rightarrow +\infty$. For that purpose, we divide the study in two cases~:
\begin{itemize}
\item Intermediate quantiles, i.e we suppose $n (1-\alpha_n) \rightarrow +\infty$. It entails that the estimation of the $\alpha_n-$quantile leads to an interpolation of sample results.
\item High quantiles. According to~\cite{dehaan1993}, we suppose $n (1-\alpha_n) \rightarrow 0$, i.e we need to extrapolate sample results to areas where no data are observed.
\end{itemize}
In both cases, the asymptotic results require some conditions we will provide throughout the section. The first one brings together the assumptions of Proposition~\ref{proploijointe}. 
\begin{itemize}
\item $(C)$ : Kernel conditions $(K1)-(K2)$ hold. In addition, $k_n \rightarrow+\infty$, $h_n \rightarrow 0$, $k_n=o(nh_n)$, $\sqrt{k_n}h_n^2 \rightarrow 0$ and $\sqrt{k_n}A \left( \frac{n}{k_n} \right) \rightarrow 0$ as $n \rightarrow +\infty$.
\end{itemize}
Condition $(C)$ will be common to both approaches, and ensures in a first time that Hill estimator is unbiased, according to Equation~\eqref{eqnormgammahat}. Moreover, $k_n=o(nh_n)$ means that $\hat{g}_{h_n}$ converges to $c_N g_N \left( M(\vec{x}) \right)$ faster than $\hat{\gamma}_{k_n}$ to $\gamma$. In practice, this condition seems appropriate, because $k_n$ must not be too large for the Hill estimator to be unbiased, and $h_n$ must be tall enough to provide a good estimation of $\ell(\vec{x})$.

\subsection{Intermediate quantiles}

We consider the case where $n (1-\alpha_n) \rightarrow +\infty$ with $\alpha_n \rightarrow 1$ as $n \rightarrow +\infty$. We recall $q_{\alpha_n}\left( Y|\vec{X}=\vec{x} \right)=\mu_{Y|\vec{X}}+\sigma_{Y|\vec{X}} \Phi_{R^*}^{-1}(\alpha_n)$. According to Equation~\eqref{equivalencephirstar}, we can approximate $\Phi_{R^*}^{-1}(\alpha_n)$ by $\Phi_R^{-1} \left(1- \left( 2+\ell(\vec{x}) \left( (1-\alpha_n)^{-1}-2 \right) \right)^{-1} \right)$. The idea is then to estimate a quantile of level $1-v_n=1- \left( 2+\ell(\vec{x}) \left( (1-\alpha_n)^{-1}-2 \right) \right)^{-1}$ on the unconditional radius $R$, easier to deal with. By noticing that $nv_n \sim \ell(\vec{x})^{-1} n(1-\alpha_n) \rightarrow +\infty $ as $n \rightarrow +\infty$, we introduce the following statistic order based estimator $\hat{q}_{\alpha_n}\left(Y|\vec{X}=\vec{x} \right)$ for $q_{\alpha_n}\left( Y|\vec{X}=\vec{x} \right)$, inspired by Theorem 2.4.1 in~\cite{dehaan2006}.

\begin{definition}[Intermediate quantile estimator]
We define $\left(\hat{q}_{\alpha_n}\left( Y|\vec{X}=\vec{x} \right) \right)_{n \in \mathbb{N}}$ as~:
\begin{equation}
\hat{q}_{\alpha_n} \left(Y|\vec{X}=\vec{x}\right)= \mu_{Y|\vec{X}}+\sigma_{Y|\vec{X}} \left( W_{\left[n \tilde{v}_n +1 \right]} \right)^{\frac{1}{\hat{\eta}_{k_n}}},
\label{eqqhat}
\end{equation}
where $\tilde{v}_n= \left( 2+\hat{\ell}_{k_n,h_n}(\vec{x}) \left( \frac{1}{1-\alpha_n}-2 \right) \right)^{-1} $, $\hat{\eta}_{k_n}$ and $\hat{\ell}_{k_n,h_n}(\vec{x})$ are respectively given in Definitions~\ref{defetahat} and~\ref{proplhat}, and $W$ is the first (or indifferently any) component of the vector $\mat{\Lambda_X}^{-1}(\vec{X}-\vec{\mu_X})$.
\label{defqhat}
\end{definition}
In order to prove the consistency of our estimator, we need a further condition $(C_{int})$ concerning the sequences $\alpha_n$ and $k_n$, usefull in the proof.
\begin{itemize}
\item $(C_{int})$ : $n (1-\alpha_n) \rightarrow +\infty$, $\ln(1-\alpha_n)=o(\sqrt{k_n})$ and $\frac{\sqrt{k_n}}{\ln(1-\alpha_n)}=o \left( \sqrt{n(1-\alpha_n)} \right)$ as $n \rightarrow +\infty$.
\end{itemize}
Obviously, $(C_{int})$ contains $n(1-\alpha_n) \rightarrow +\infty$, as mentioned above. Furthermore, $\ln(1-\alpha_n)=o(\sqrt{k_n})$ ensures that the rate of convergence in Theorem~\ref{propnormint1} goes to infinity (see below) and the last relationship allows us to eliminate a term in the proof. In order to make this condition more meaningful, let us propose a simple example: we choose our sequences in polynomial forms $k_n=n^b$, $0<b<1$ and $\alpha_n=1-n^{-a}, a>0$. It is straightforward to see that $\ln(1-\alpha_n)=o(k_n)$ and $\ln(n(1-\alpha_n))=o(k_n)$ $, \forall a>0, 0<b<1$. However, $\frac{\sqrt{k_n}}{\ln(1-\alpha_n)}=o \left( \sqrt{n(1-\alpha_n)} \right)$ if and only if $a<1$, i.e $n(1-\alpha_n) \rightarrow +\infty$ as $n \rightarrow +\infty$. \\ In a first time, we give a result concerning the asymptotic behavior of $\hat{q}_{\alpha_n} \left(Y|\vec{X}=\vec{x} \right)$ with respect to $q_{\alpha_n \uparrow}\left( Y|\vec{X}=\vec{x} \right)$. Then, with Equation~\eqref{equivalencephirstar}, we easily deduce a consistency result for $\hat{q}_{\alpha_n}\left( Y|\vec{X}=\vec{x} \right) $.

\begin{theoreme}[Consistency of $\hat{q}_{\alpha_n} \left(Y|\vec{X}=\vec{x}\right)$]
Let us denote $v_n=\left( 2+\ell(\vec{x}) \left( (1-\alpha_n)^{-1}-2 \right) \right)^{-1}$ and $\tilde{v}_n=\left( 2+\hat{\ell}_{k_n,h_n}(\vec{x}) \left( (1-\alpha_n)^{-1}-2 \right) \right)^{-1}$. Under Assumption~\ref{hyp2}, and conditions $(C),(C_{int})$~: 
\begin{equation}
\frac{\sqrt{k_n}}{\ln \left( 1-\alpha_n \right)} \left( \frac{\hat{q}_{\alpha_n} \left( Y|\vec{X}=\vec{x} \right)}{q_{\alpha_n \uparrow}\left( Y|\vec{X}=\vec{x} \right)}-1 \right) \underset{n \rightarrow +\infty}{\rightarrow} \mathcal{N} \left( 0,\frac{N^2 \gamma^4}{\left( \gamma N+1 \right)^4} \right).
\label{eqnormint1}
\end{equation}
And therefore~:
\begin{equation}
\frac{\hat{q}_{\alpha_n} \left( Y|\vec{X}=\vec{x} \right)}{q_{\alpha_n}\left( Y|\vec{X}=\vec{x} \right)} \overset{\mathbb{P}}{\rightarrow} 1.
\label{eqcvprobainterm}
\end{equation}
\label{propnormint1}
\end{theoreme}

The same asymptotic normality with $\Phi_{R^*}^{-1}(\alpha_n)$ instead of $\Phi_R^{-1} \left( 1-v_n \right)^{\frac{1}{\eta}}$ may be deduced from Proposition~\ref{propnormint1} under the condition $$\underset{n \rightarrow +\infty}{\lim} \text{ } \frac{\sqrt{k_n}}{\ln \left( 1-\alpha_n \right)} \ln \left( \frac{\Phi_{R}^{-1} \left( 1-v_n \right)^{\frac{1}{\eta}}}{\Phi_{R^*}^{-1}(\alpha_n)} \right) = 0.$$ This condition, which seems quite simple, is difficult to prove in a general context. Indeed, we need a second order expansion of Equation~\eqref{equivalencephirstar}. But the second order properties of the unconditional quantile $\Phi_R^{-1}$ given by Assumption~\ref{hyp2} are not necessarily the same as those of the conditional quantile $\Phi_{R^*}^{-1}$, which makes the study complicated. However, in some simple cases, we are able to solve the problem. We thus give another assumption, stronger that Assumption~\ref{hyp2}. In the following, we refer to this assumption for results of asymptotic normality.

\begin{assumption}
$\forall d \in \mathbb{N}^*$, there exists $\lambda_1,\lambda_2 \in \mathbb{R}$ such that~:
\begin{equation}
c_d g_d(t)=\lambda_1 t^{-\frac{d+\gamma^{-1}}{2}} \left[ 1+\lambda_2 t^{\frac{\rho}{2  \gamma}}+o\left( t^{\frac{\rho}{2  \gamma}} \right) \right].
\label{eqhypforte}
\end{equation}
\label{hypforte}
\end{assumption}
It is obvious that Assumption~\ref{hypforte} implies Assumption~\ref{hyp2}. Indeed, according to~\cite{hua}, Equation~\eqref{eqhypforte} is equivalent to say that $c_1 g_1(t^2)$ is regularly varying of second order with indices $-1-\gamma^{-1}$, $\rho/\gamma$ and an auxiliary function proportional to $t^{\frac{\rho}{\gamma}}$. Then, Proposition 6 in~\cite{hua} entails $\bar{\Phi}_R(t)$ is second order regularly varying with $-\gamma^{-1}$, $\rho/\gamma$ and the same kind of auxiliary function. Finally, this is equivalent (see~\cite{dehaan2006}) to Assumption~\ref{hyp2} with indicated $\gamma$ and $\rho$, and an auxiliary function $A(t)$ proportional to $t^{\rho}$. \\ Furthermore, according to~\cite{kano}, the dependance on $d$ in Equation~\eqref{eqhypforte} remains coherent with the assumption of consistent elliptical distributions, the latter having to have a function $g_d$ depending on $d$. As an example, the Student distribution fills Assumption~\ref{hypforte}. The latter allows us to provide a second order expansion for Equation~\eqref{equivalencephirstar}. \\ In order to prove the asymptotic normality of $\hat{q}_{\alpha_n} \left(Y|\vec{X}=\vec{x}\right)$, we add a technical condition $\left(C_{int}^{HG} \right)$ that involves tail indices $\gamma$ and $\rho$.
\begin{itemize}
\item $\left(C_{int}^{HG} \right)$ : $(C_{int})$ holds. In addition, $\sqrt{k_n}(1-\alpha_n)=o \left( \ln(1-\alpha_n) \right)$, and~:
\begin{equation}
\underset{n \rightarrow +\infty}{\lim} \text{ } \frac{\sqrt{k_n}}{\ln \left( 1-\alpha_n \right)} \left( 1-\alpha_n \right)^{\frac{\min(-\rho,2 \gamma)}{\gamma N+1}}=0.
\label{conditioncint}
\end{equation}
\end{itemize}
Condition $\left(C_{int}^{HG} \right)$ means that sequence $k_n$ must not be too large. In view of Equation~\eqref{conditioncint}, it is obvious that if $N$ or $\gamma$ goes to infinity, $\left(C_{int}^{HG} \right)$ is not filled. The tail of the underlying distribution may thus not be too heavy, and the size $N$ of the covariate not too large. Similarly, they no longer hold if $\gamma$ or $\rho$ goes to 0, i.e. if the underlying distribution is either too lightly varying, or its c.d.f. takes too long to behave like $\lambda t^{-1/\gamma}$.

\begin{prop}[Asymptotic normality of $\hat{q}_{\alpha_n} \left(Y|\vec{X}=\vec{x}\right)$]
Assume that Assumption~\ref{hypforte} and conditions $(C), \left( C_{int}^{HG} \right)$ hold. Then~:
\begin{equation}
\frac{\sqrt{k_n}}{\ln \left( 1-\alpha_n \right)} \left( \frac{\hat{q}_{\alpha_n} \left(Y|\vec{X}=\vec{x}\right)}{q_{\alpha_n} \left(Y|\vec{X}=\vec{x}\right)}-1 \right) \underset{n \rightarrow +\infty}{\rightarrow} \mathcal{N} \left( 0,\frac{N^2 \gamma^4}{\left( \gamma N+1 \right)^4} \right).
\label{eqnormint2}
\end{equation}
\label{propnormint2}
\end{prop}

We notice that asymptotic variance in Equation~\eqref{eqnormint1} tends to 0 as the number of covariates $N$ goes to $+\infty$. Indeed, we observe a fast convergence of $\hat{q}_{\alpha_n}$ to $q_{\alpha_n \uparrow}$ when $N$ is large. However, $\left(  C_{int}^{HG} \right)$ is not filled if $N$ is tall. Then asymptotic normality~\eqref{eqnormint2} no longer holds. This is explained by the fact that more $N$ is tall, more $q_{\alpha_n}\left( Y|\vec{X}=\vec{x} \right)/q_{\alpha_n \uparrow}\left( Y|\vec{X}=\vec{x} \right)$ (see Equation~\eqref{equivalencephirstar}) tends to 1 slowly.

\subsection{High quantiles}

We now consider $n (1-\alpha_n) \rightarrow 0$ as $n \rightarrow +\infty$. In the following definition, we introduce another quantile estimator $\hat{\hat{q}}_{\alpha_n} \left(Y|\vec{X}=\vec{x} \right)$ for $q_{\alpha_n} \left(Y|\vec{X}=\vec{x} \right)$. 
We first recall that the idea is to estimate an unconditional quantile of level $1-v_n=1-\left( 2+\ell(\vec{x}) \left( (1-\alpha_n)^{-1}-2 \right) \right)^{-1}$. A quick calculation proves that $v_n$ is asymptotically equivalent to $\ell(\vec{x})^{-1}(1-\alpha_n)$, and therefore $nv_n \rightarrow 0$ as $n \rightarrow +\infty$. The use of statistic order (at level $nv_n$) is then impossible in that case. According to Theorem 4.3.8 in~\cite{dehaan2006}, a way to estimate such a quantile may be to take the statistic order at the intermediate level $k_n$ (we recall $k_n \rightarrow +\infty$), and apply an extrapolation coefficient $\left(k_n/(n v_n) \right)^{\gamma}$. This approach inspired the following estimator.  
\begin{definition}[High quantile estimator]
We define $\left(\hat{\hat{q}}_{\alpha_n}\left( Y|\vec{X}=\vec{x} \right) \right)_{n \in \mathbb{N}}$ as~:
\begin{equation}
\hat{\hat{q}}_{\alpha_n} \left(Y|\vec{X}=\vec{x} \right)=\mu_{Y|\vec{X}}+\sigma_{Y|\vec{X}} \left[W_{[k_n+1]} \left( \frac{k_n}{n} \left( 2+\hat{\ell}_{k_n,h_n}(\vec{x}) \left( \frac{1}{1-\alpha_n}-2 \right) \right) \right)^{\hat{\gamma}_{k_n}} \right]^{\frac{1}{\hat{\eta}_{k_n}}}.
\label{eqqhathat}
\end{equation}
\label{defqhathat}
\end{definition}
The aim is now to study the asymptotic properties of $\hat{\hat{q}}_{\alpha_n} \left( Y|\vec{X}=\vec{x} \right)$. As for the intermediate quantile estimator, we propose a result of asymptotic normality, under a condition $(C_{high})$ (given below) which we then refine under Assumption~\ref{hypforte}. 
\begin{itemize}
\item $(C_{high})$ : $n (1-\alpha_n) \rightarrow 0$, $\ln \left(n(1-\alpha_n) \right)=o(\sqrt{k_n})$ and $ \frac{\ln \left( 1-\alpha_n \right)}{\ln \left( \frac{n}{k_n}(1-\alpha_n) \right)} \rightarrow \theta \in [0,+\infty[$ as $n \rightarrow +\infty$.
\end{itemize}
The second statement is added in order to apply Theorem 4.3.8 in~\cite{dehaan2006}, and the third one is a notation used in the following. Let us propose a simple example: if we choose our sequences in polynomial forms $k_n=n^b$, $0<b<1$ and $\alpha_n=1-n^{-a}, a>0$, the first condition is filled if and only if $a>1$, $\ln(n(1-\alpha_n))=o(\sqrt{k_n})$ and the last assertion holds with a particular $\theta$ given later. \\ The consistency result that follows immediatly is given just below.

\begin{theoreme}[Consistency of high quantile estimator]
Let us denote $v_n=\left( 2+\ell(\vec{x}) \left( (1-\alpha_n)^{-1}-2 \right) \right)^{-1}$ and $\tilde{v}_n=\left( 2+\hat{\ell}_{k_n,h_n}(\vec{x}) \left( (1-\alpha_n)^{-1}-2 \right) \right)^{-1}$. Under Assumption~\ref{hyp2}, and conditions $(C),(C_{high})$~:
\begin{equation}
\frac{\sqrt{k_n}}{\ln \left( \frac{k_n}{n (1-\alpha_n)} \right)} \left( \frac{\hat{\hat{q}}_{\alpha_n} \left(Y|\vec{X}=\vec{x}  \right)}{q_{\alpha_n \uparrow} \left(Y|\vec{X}=\vec{x}  \right)}-1 \right) \underset{n \rightarrow +\infty}{\rightarrow} \mathcal{N} \left( 0, \left( \frac{\gamma}{\gamma N+1}-\theta \frac{N \gamma^2}{(\gamma N+1)^2} \right)^2 \right).
\label{eqasympnormphirstar}
\end{equation}
And therefore~:
\begin{equation}
\frac{\hat{\hat{q}}_{\alpha_n} \left( Y| \vec{X}=\vec{x} \right)}{q_{\alpha_n} \left( Y| \vec{X}=\vec{x} \right)} \overset{\mathbb{P}}{\rightarrow} 1 \text{ as } n \rightarrow +\infty.
\end{equation}
\label{thasympnormphirstar}
\end{theoreme}

We can emphasize that condition $(C_{high})$ is filled in most of the common cases. Indeed, the simple examples to find that do not satisfy (ii) are of the form $\alpha_n=1-n^{-1} \ln(n)^{-\kappa}, \kappa >0$ and $k_n=\ln(n)$. But such a choice of sequences would lead to a poor estimation of $\hat{\gamma}_{k_n}$ and $\hat{\eta}_{k_n}$, since $k_n \rightarrow +\infty$ very slowly, and moreover a poor estimation of the quantile, the level $\alpha_n$ tending to 1 slowly. These sequences are therefore not recommanded in practice. Next corollary gives the value of $\theta$ when sequences $k_n$ and $\alpha_n$ have a polynomial form.

\begin{cor}
Under Assumption~\ref{hyp2}, conditions $(C), (C_{high})$, and taking $k_n=n^b, 0<b<1$ and $\alpha_n=1-n^{-a},a>1$, asymptotic relationship~\eqref{eqasympnormphirstar} holds with $\theta=\frac{a}{a+b-1}$.
\label{corollaire}
\end{cor}

As for the intermediate quantile estimator, asymptotic normality~\eqref{eqasympnormphirstar} may be improved under the condition $$\underset{n \rightarrow +\infty}{\lim} \text{ } \frac{\sqrt{k_n}}{\ln \left( \frac{k_n}{n(1-\alpha_n)} \right)} \ln \left( \frac{\Phi_{R}^{-1} \left( 1-v_n \right)^{\frac{1}{\eta}}}{\Phi_{R^*}^{-1}(\alpha_n)} \right) = 0$$ Assumption~\ref{hypforte} places us in a framework where it is quite simple to prove it, if we add the following condition~:
\begin{itemize}
\item $\left(C_{high}^{HG} \right)$ : $(C_{high})$ holds. In addition,
\begin{equation}
\underset{n \rightarrow +\infty}{\lim} \text{ } \frac{\sqrt{k_n}}{\ln \left( \frac{k_n}{n (1-\alpha_n)} \right)} \left( 1-\alpha_n \right)^{\frac{\min(-\rho,2\gamma)}{\gamma N+1}}=0.
\label{conditionchigh}
\end{equation}
\end{itemize}
As $\left(C_{int}^{HG} \right)$, condition $\left(C_{high}^{HG} \right)$ means that sequence $k_n$ must be small enough. In view of Equation~\eqref{conditionchigh}, we deduce that if $N$ or $\gamma$ goes to infinity, $\left(C_{high}^{HG} \right)$ is not filled. The tail of the underlying distribution may thus not be too heavy, and the size $N$ of the covariate not too large. Similarly, they no longer hold if $\gamma$ or $\rho$ goes to 0. \\
By combining Assumption~\ref{hypforte} and $\left(C_{high}^{HG} \right)$, the following result is obtained.

\begin{prop}[Asymptotic normality of high quantile estimator]
Assume that Assumption~\ref{hypforte} and conditions $(C), \left(C_{high}^{HG} \right)$ hold. Then~:
\begin{equation}
\frac{\sqrt{k_n}}{\ln \left( \frac{k_n}{n (1-\alpha_n)} \right)} \left( \frac{\hat{\hat{q}}_{\alpha_n} \left(Y|\vec{X}=\vec{x} \right)}{q_{\alpha_n} \left(Y|\vec{X}=\vec{x} \right)}-1 \right) \underset{n \rightarrow +\infty}{\rightarrow} \mathcal{N} \left( 0,\left( \frac{\gamma}{\gamma N+1}-\theta \frac{N \gamma^2}{(\gamma N+1)^2} \right)^2 \right).
\label{eqasympnormphirstar2}
\end{equation}
\label{propasympnormphirstar2}
\end{prop}

We can make the same kind of remark as in the previous subsection when $N$ is large. In the following, we give estimators for two other classes of extreme risk measures, based on the estimators given in Equations~\eqref{eqqhat} and \eqref{eqqhathat}. The first one generalizes quantiles.

\section{Some extreme risk measures estimators} \label{lpquantiles}

\subsection{$L_p-$quantiles}

Let $Z$ be a real random variable. The $L_p-$quantiles of $Z$ with level $\alpha \in ]0,1[$ and $p>0$, denoted $q_{p,\alpha}(Z)$, is solution of the minimization problem (see~\cite{chen})~:
\begin{equation}
q_{p,\alpha}(Z)=  \underset{z \in \mathbb{R}}{\arg \min} \text{ } \mathbb{E} \left[ (1-\alpha) \left( z-Z \right)_+^p + \alpha \left( Z-z \right)_+^p \right],
\label{eqlpquantiles}
\end{equation}
where $Z_+=Z \mathds{1}_{ \{Z>0 \}}$. According to \cite{koenker}, the case $p=1$ leads to the quantile $q_{1,\alpha}(Z)=F_Z^{-1}(\alpha)$, where $F_Z$ is the c.d.f of $Z$. The case $p=2$, formalized in~\cite{newey}, leads to more complicated calculations, and admits, with the exception of some particular cases (see, e.g., \cite{koenkerexpectile1}), no general formula. The general case $p \geq 1$ has seen some recent advances. \cite{belliniklar} has shown that $L_p-$quantiles get the translation equivariance and positively homogeneity properties for $p>1$. More recently, the particular case of Student distributions has, for example, been explored in~\cite{bignozzi}. However, it seems difficult to obtain a general formula. On the other hand, in the case of extreme levels $\alpha$, i.e. when $\alpha$ tends to 1,~\cite{girardlp} proved that the following relationship holds, for a heavy-tailed random variable with tail index $\gamma$.
\begin{equation}
\frac{q_{p,\alpha} \left(Z \right)}{q_{\alpha}(Z)} \underset{\alpha \rightarrow 1}{\rightarrow} \left[ \frac{\gamma}{B \left( p, \gamma^{-1}-p+1 \right)} \right]^{-\gamma} := f_L \left( \gamma,p \right),
\label{equivalencelp}
\end{equation}  
where $B(.,.)$ is the beta function. We add that for a Pareto-type distribution with tail index $\gamma$, the $L_p-$quantile exists if and the only if the moment of order $p-1$ exists, i.e. if $\gamma<1/p$. The expectile case $p=2$ leads to the result of~\cite{belliniklar}. Using this result, we can estimate the conditional $L_p-$quantiles from the quantile estimated in Section~\ref{quantiles}. For that purpose, we need to know the tail index of the conditional radius $R^*$, given in the following lemma.
\begin{lemme}
The conditional distribution $Y|\vec{X}=\vec{x}$ is attracted to a maximum domain of Pareto-type distribution with tail index $(\gamma^{-1}+N)^{-1}$, i.e 
\begin{equation}
\underset{t \rightarrow +\infty}{\lim} \frac{\bar{\Phi}_{R^*}(\omega t)}{\bar{\Phi}_{R^*}(t)}=\omega^{-\frac{1}{\gamma}-N}.
\end{equation}
\label{lemmegammacond}
\end{lemme}

With Lemma~\ref{lemmegammacond} and Equation~\eqref{equivalencelp}, we define the following estimators for the $L_p-$quantile of $Y|\vec{X}=\vec{x}$, according to whether if $n(1-\alpha_n)$ tends to $0$ or $+\infty$.

\begin{definition}
Let $(\alpha_n)_{n \in \mathbb{N}}$ be a sequence such that $\alpha_n \rightarrow 1$ as $n \rightarrow +\infty$. If either $p \leq N$ or $\gamma<\frac{1}{p-N}$, we define:
\begin{equation}
\left \{
\begin{array}{cc}
      \hat{q}_{p,\alpha_n} \left(Y|\vec{X}=\vec{x} \right) =& \mu_{Y|\vec{X}}+\sigma_{Y|\vec{X}} \left( W_{\left[n \tilde{v}_n +1 \right]} \right)^{1/\hat{\eta}_{k_n}} f_L \left( \left( \hat{\gamma}_{k_n}^{-1}+N \right)^{-1}, p \right) \\
      \hat{\hat{q}}_{p,\alpha_n} \left(Y|\vec{X}=\vec{x} \right) =& \mu_{Y|\vec{X}}+\sigma_{Y|\vec{X}} \left[W_{[k_n+1]} \left( \frac{k_n}{n \tilde{v}_n}  \right)^{\hat{\gamma}_{k_n}} \right]^{1/\hat{\eta}_{k_n}} f_L \left( \left( \hat{\gamma}_{k_n}^{-1}+N \right)^{-1}, p \right)  \\
\end{array}.
\right.
\label{eqqlphat}
\end{equation}
where $\hat{\gamma}_{k_n}$ and $\tilde{v}_n$ are respectively given in Equation~\eqref{hillest} and Theorem~\ref{propnormint1}.
\label{deflphat}
\end{definition}
We have proved the convergence in probability of $\hat{q}_{\alpha_n}\left(Y|\vec{X}=\vec{x} \right)$ and $\hat{\hat{q}}_{\alpha_n}\left(Y|\vec{X}=\vec{x} \right)$. Furthermore, the convergence in probability of the asymptotic term, and consequently the empirical $L_p-$quantile is not difficult to get, this is why we omit the proof.

\begin{prop}[Consistency of $L_p-$quantile estimators]
Assume that Assumption~\ref{hyp2} and condition $(C)$ hold. Under conditions $(C_{int})$ and $(C_{high})$ respectively, $\hat{q}_{p,\alpha_n}\left(Y|\vec{X}=\vec{x} \right)$ and $\hat{\hat{q}}_{p,\alpha_n}\left(Y|\vec{X}=\vec{x} \right)$ are consistent, i.e.~:
\begin{equation}
\left \{
\begin{array}{cc}
       \frac{\hat{q}_{p,\alpha_n} \left(Y|\vec{X}=\vec{x} \right)}{q_{p,\alpha_n}\left(Y|\vec{X}=\vec{x} \right)} \overset{\mathbb{P}}{\rightarrow}& 1  \\
        \frac{\hat{\hat{q}}_{p,\alpha_n}\left(Y|\vec{X}=\vec{x} \right)}{q_{p,\alpha_n}\left(Y|\vec{X}=\vec{x} \right)} \overset{\mathbb{P}}{\rightarrow}& 1  \\
\end{array}.
\right.
\label{eqconsistencelp}
\end{equation}
\end{prop}

Using the second order expansion of Equation~\eqref{equivalencelp} given in~\cite{girardlp}, and doing some stronger assumptions, we can deduce the following asymptotic normality results. For that purpose, let us add two conditions.
\begin{itemize}
\item $\left(C_{int}^{L_p}\right)$ : $(C_{int})$ holds. In addition, $\sqrt{k_n}(1-\alpha_n) =o \left( \ln(1-\alpha_n) \right)$, and~:
\begin{equation}
\underset{n \rightarrow +\infty}{\lim} \text{ } \frac{\sqrt{k_n}}{\ln \left( 1-\alpha_n \right)} \left( 1-\alpha_n \right)^{\frac{\min(-\rho,\gamma)}{\gamma N+1}}=0.
\end{equation}
\item $\left(C_{high}^{L_p}\right)$ : $(C_{high})$ holds. In addition,
\begin{equation}
\underset{n \rightarrow +\infty}{\lim} \text{ } \frac{\sqrt{k_n}}{\ln \left( \frac{k_n}{n (1-\alpha_n)} \right)} \left( 1-\alpha_n \right)^{\frac{\min(-\rho,\gamma)}{\gamma N+1}}=0.
\end{equation}
\end{itemize} 
These conditions will be used below. If we compare $\left(C_{int}^{L_p} \right)$ and $\left(C_{int}^{L_p} \right)$ with $\left(C_{int}^{HG} \right)$ and $\left(C_{int}^{HG} \right)$ respectively, sequence $k_n$ must be chosen smaller. Finally, we can draw the same conclusions than above, i.e. these conditions are applicable for regularly varying distributions with an intermediate level $\gamma$, and a small number of covariates $N$. \\
To sum up, among all these conditions, we can deduce the following ordering~:$$ \left \{
\begin{array}{cc}
      \left(C_{int}^{L_p}\right) \Rightarrow \left(C_{int}^{HG}\right) \Rightarrow \left(C_{int} \right) \\
       \left(C_{high}^{L_p}\right) \Rightarrow \left(C_{high}^{HG}\right) \Rightarrow \left(C_{high} \right) \\
\end{array}.
\right. $$

\begin{prop}[Asymptotic normality of $L_p-$quantile estimators]
Assume that Assumption~\ref{hypforte} and condition $(C)$ hold. Under conditions $\left(C_{int}^{L_p}\right)$ and $\left(C_{high}^{L_p}\right)$ respectively, and if $p>1$, then~:
\begin{equation}
\left \{
\begin{array}{cc}
        \frac{\sqrt{k_n}}{\ln \left( 1-\alpha_n \right)} \left( \frac{\hat{q}_{p,\alpha_n}\left( Y|\vec{X}=\vec{x} \right)}{q_{p,\alpha_n}\left( Y|\vec{X}=\vec{x} \right)}-1 \right) \underset{n \rightarrow +\infty}{\rightarrow} & \mathcal{N} \left( 0, \frac{N^2 \gamma^4}{\left( \gamma N+1 \right)^4} \right) \\
         \frac{\sqrt{k_n}}{\ln \left( \frac{k_n}{n (1-\alpha_n)} \right)} \left( \frac{\hat{\hat{q}}_{p,\alpha_n}\left( Y|\vec{X}=\vec{x} \right)}{q_{p,\alpha_n}\left( Y|\vec{X}=\vec{x} \right)}-1 \right) \underset{n \rightarrow +\infty}{\rightarrow} & \mathcal{N} \left( 0,\left( \frac{\gamma}{\gamma N+1}-\theta \frac{N \gamma^2}{(\gamma N+1)^2} \right)^2\right) \\
\end{array}.
\right.
\label{eqnormlp}
\end{equation}
\label{propnormlp}
\end{prop}

An example of $L_2-$quantile, or expectile, is provided in Section~\ref{simstudy}. The second risk measure we focus on is called Haezendonck-Goovaerts risk measure.

\subsection{Haezendonck-Goovaerts risk measures} Let $Z$ be a real random variable, and $\varphi$ a non negative and convex function with $\varphi(0)=0$, $\varphi(1)=1$ and $\varphi(+\infty)=+\infty$. The Haezendonck-Goovaerts risk measure of $Z$ with level $\alpha \in ]0,1[$ associated to $\varphi$, is given by the following (see~\cite{tangyang})~:
\begin{equation}
H_{\alpha}(Z)= \underset{z \in \mathbb{R}}{\inf} \left\{ z+H_{\alpha}(Z,z) \right\},
\end{equation}
where $H_{\alpha}(Z,z)$ is the unique solution $h$ to the equation~:
\begin{equation}
\mathbb{E} \left[ \varphi \left( \frac{(Z-z)_+}{h} \right) \right]=1-\alpha.
\end{equation}
$\varphi$ is called Young function. This family of risk measures has been firstly introduced as Orlicz risk measure in~\cite{haezendonck}, then Haezendonck risk measure in~\cite{goovaerts}, and finally Haezendonck-Goovaerts risk measure in~\cite{tangyang}. According to~\cite{bellinihaezendonck}, such a risk measure is coherent, and therefore translation equivariant and positively homogeneous. The particular case $\varphi(t)=t$ leads to the Tail Value at Risk with level $\alpha$ TVaR$_{\alpha}(X)$, introduced in~~\cite{artzner}. In the following, we denote $H_{p,\alpha}(Z)$ the Haezendonck-Goovaerts risk measure of $Z$ with a power Young function $t^p, p \geq 1$. In~\cite{tangyang}, the authors provided the following result.
\begin{prop}[\cite{tangyang}]
If $Z$ fills Assumption~\ref{hyp2}, and taking a Young function $\varphi(t)=t^p, p\geq 1$, then the following relationship holds~:
\begin{equation}
\frac{H_{p,\alpha}(Z)}{q_{\alpha}(Z)} \underset{\alpha \rightarrow 1}{\rightarrow} \frac{\gamma^{-1} \left(\gamma^{-1}-p \right)^{p \gamma-1}}{p^{\gamma(p-1)}} B \left( \gamma^{-1}-p,p \right)^{\gamma} := f_H \left( \gamma, p \right).
\label{eqequivalencehaezendonck}
\end{equation}
\label{propequivalencehaezendonck}
\end{prop}
In particular, taking $p=1$ leads to TVaR$_{\alpha}(Z) \sim (1-\gamma)^{-1}$ $q_{\alpha}(Z)$ as $\alpha \rightarrow 1$. Using Lemma~\ref{lemmegammacond}, extreme quantiles estimators in Definitions~\ref{defqhat}, \ref{defqhathat} and Proposition~\ref{propequivalencehaezendonck}, we can deduce estimators for extreme Haezendonck-Goovaerts risk measure $H_{p,\alpha} \left( Y| \vec{X}=\vec{x} \right)$ (with power Young function $\varphi(t)=t^p, p \geq 1$) of $Y| \vec{X}=\vec{x}$.
\begin{definition}
Let $(\alpha_n)_{n \in \mathbb{N}}$ be a sequence such that $\alpha_n \rightarrow 1$ as $n \rightarrow +\infty$. If either $p \leq N$ or $\gamma<\frac{1}{p-N}$, we define~:
\begin{equation}
\left \{
\begin{array}{cc}
      \hat{H}_{p,\alpha_n}\left(Y| \vec{X}=\vec{x} \right) =& \mu_{Y|\vec{X}}+\sigma_{Y|\vec{X}} \left( W_{\left[n \tilde{v}_n +1 \right]} \right)^{1/\hat{\eta}_{k_n}} f_H \left( \left( \hat{\gamma}_{k_n}^{-1}+N \right)^{-1}, p \right) \\
      \hat{\hat{H}}_{p,\alpha_n}\left(Y| \vec{X}=\vec{x} \right) =& \mu_{Y|\vec{X}}+\sigma_{Y|\vec{X}} \left[W_{[k_n+1]} \left( \frac{k_n}{n \tilde{v}_n}  \right)^{\hat{\gamma}_{k_n}} \right]^{1/\hat{\eta}_{k_n}} f_H \left( \left( \hat{\gamma}_{k_n}^{-1}+N \right)^{-1}, p \right)  \\
\end{array}.
\right.
\label{eqhhat}
\end{equation}
\label{defhhat}
\end{definition}
The condition $p \leq N$ or $\gamma<\frac{1}{p-N}$ simply ensures the existence of $H_{p,\alpha_n}\left(Y| \vec{X}=\vec{x} \right)$. Using the consistency results given in Propositions~\ref{propnormint1} and~\ref{thasympnormphirstar}, the consistency of these estimators is immediate. The proof is also omitted from the appendix.
\begin{prop}[Consistency of H-G estimators]
Assume that Assumption~\ref{hyp2} and condition $(C)$ hold. Under conditions $(C_{int})$ and $(C_{high})$ respectively, $\hat{H}_{p,\alpha_n}\left(Y| \vec{X}=\vec{x} \right)$ and $\hat{\hat{H}}_{p,\alpha_n}\left(Y| \vec{X}=\vec{x} \right)$ are consistent, i.e.~:
\begin{equation}
\left \{
\begin{array}{cc}
       \frac{\hat{H}_{p,\alpha_n}\left(Y| \vec{X}=\vec{x} \right)}{H_{p,\alpha_n}\left(Y| \vec{X}=\vec{x} \right)} \overset{\mathbb{P}}{\rightarrow}& 1  \\
        \frac{\hat{\hat{H}}_{p,\alpha_n}\left(Y| \vec{X}=\vec{x} \right)}{H_{p,\alpha_n}\left(Y| \vec{X}=\vec{x} \right)} \overset{\mathbb{P}}{\rightarrow}& 1  \\
\end{array}.
\right.
\label{eqconsistencehaezendonck}
\end{equation}
\end{prop}

\begin{prop}[Asymptotic normality of H-G estimators]
Assume that Assumption~\ref{hypforte} and condition $(C)$ hold. Under conditions $\left( C_{int}^{HG} \right)$ and $\left( C_{high}^{HG} \right)$ respectively, we have~:
\begin{equation}
\left \{
\begin{array}{cc}
        \frac{\sqrt{k_n}}{\ln \left( 1-\alpha_n \right)} \left( \frac{\hat{H}_{p,\alpha_n}\left(Y| \vec{X}=\vec{x} \right)}{H_{p,\alpha_n}\left(Y| \vec{X}=\vec{x} \right)}-1 \right) \underset{n \rightarrow +\infty}{\rightarrow} & \mathcal{N} \left( 0, \frac{N^2 \gamma^4}{\left( \gamma N+1 \right)^4} \right) \\
         \frac{\sqrt{k_n}}{\ln \left( \frac{k_n}{n (1-\alpha_n)} \right)} \left( \frac{\hat{\hat{H}}_{p,\alpha_n}\left(Y| \vec{X}=\vec{x} \right)}{H_{p,\alpha_n}\left(Y| \vec{X}=\vec{x} \right)}-1 \right) \underset{n \rightarrow +\infty}{\rightarrow} & \mathcal{N} \left( 0,\left( \frac{\gamma}{\gamma N+1}-\theta \frac{N \gamma^2}{(\gamma N+1)^2} \right)^2 \right) \\
\end{array}.
\right.
\label{eqnormhg}
\end{equation}
\label{propnormhg}
\end{prop}

We can emphasize that conditions for asymptotic normality are less strong in the case of Haezendonck-Goovaerts risk measures. We propose some examples (with $p=1$, i.e TVaR) in Sections~\ref{simstudy} and~\ref{real}.

\section{Simulation study} \label{simstudy}

In this section, we apply our estimators to 100 samples of $n$ simulations of a Student vector $\vec{Z}=(\vec{X},Y) \in \mathbb{R}^4$ ($\vec{X} \in \mathbb{R}^3$ and $Y \in \mathbb{R}$) with $\nu=2$ degrees of freedom, and compare with theoretical results. According to~\cite{dehaan2006}, the Student distribution with $\nu$ degrees of freedom fills Assumption~\ref{hyp2} with indices $\gamma=1/\nu$, $\rho=-2/\nu$, and an auxiliary function $A(t)$ proportional to $t^{-2/\nu}$. The latter even fills Assumption~\ref{hypforte}, and is the only heavy-tailed elliptical distribution (to our knowledge) where we can obtain closed formula for conditional quantiles. In addition, such a degree of freedom makes the tail of the distribution sufficiently heavy to easily observe the asymptotic results. We can notice that the unconditional distribution $Y$ has tail index $1/2$, then, using Lemma~\ref{lemmegammacond}, the conditional distribution $Y|\vec{X}=\vec{x}$ has tail index $2/7<1/2$, and admits quantile, expectile ($L_2-$quantile) and TVaR. This section beeing uniquely devoted to the performance of our estimators, we take for conveniance $\vec{\mu}=\vec{0}_{\mathbb{R}^4}$ and $\mat{\Sigma}=\mat{I}_4$. Let us now estimate the extreme quantiles of $Y|\vec{X}=\vec{x}$. For that purpose, we have to chose an arbitrary value of $\vec{x}$. We thus suppose for example that the observed covariates $\vec{x}$ satisfy $M(\vec{x})=1$.

\subsection{Choice of parameters}

As mentioned in Sections~\ref{parametres} and~\ref{quantiles}, the asymptotic results obtained are sensitive to the choice of sequences $k_n$, $h_n$, $\alpha_n$, and to a lesser extent to the kernel $K$. The latter will be the gaussian p.d.f in the following. Concerning the sequences, we propose in this section to consider the polynomial forms $\alpha_n=1-n^{-a}$, $a>0$, $k_n=n^b$, $b>0$ and $h_n=n^{-c}$, $c>0$. In order to deal with high quantiles, we fix in a first time $a=1.25$. We now have to chose carefully the parameters $b$ and $c$, fulfilling the conditions $(C)$, $\left( C_{high} \right)$ and $\left( C_{high}^{HG} \right)$. $(C)$ imposes $b<1-c$, $b<4c$ and $b<4/(\nu+4)=2/3$, $\left( C_{high} \right)$ is satisfied with $\theta=a/(a+b-1)$ (see Corollary~\ref{corollaire}), $\left( C_{high}^{HG} \right)$ entails $b \leq 2a =2.5$ and $b \leq 4a/(N+\nu)=1$. Finally, it seems reasonable to chose $b$ (respectively $c$) as tall (respectively small) as possible. The choices $b=0.6$ and $c=0.2$ seem to be a good compromise.

\subsection{Extremal parameters estimation}

The next step is to estimate the quantities $\eta$ and $\ell(\vec{x})$. For that purpose, we use our estimators $\hat{\eta}_{k_n}$ and $\hat{\ell}_{k_n,h_n}(\vec{x})$ respectively introduced in Equations~\eqref{eqetahat} and~\eqref{eqlhat}. These two estimators are related to the Hill estimator $\hat{\gamma}_{k_n}$, and asymptotic results of Section~\ref{parametres} hold only if the data is independent. This is why we do the estimation of $\gamma$ only with the $n$ realizations of the first component from the vector $\vec{Z}$.

\begin{figure}[h!]
\begin{center}
\includegraphics[scale=0.55]{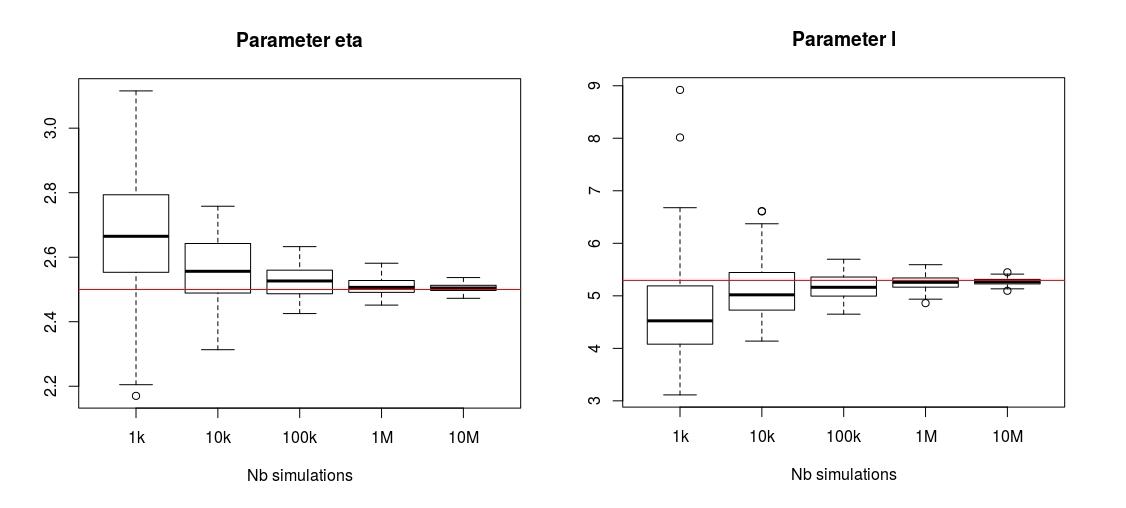}
\end{center}
\caption{From left to right: boxplots of 100 estimators $\hat{\eta}_{k_n}$ and $\hat{\ell}_{k_n,h_n}(\vec{x})$, for different sample sizes $n$. Theoretical values are in red. The chosen sequences are $k_n=n^{0.6}$, $h_n=n^{-0.2}$ and $\alpha_n=1-n^{-1.25}$.}
\label{plot_sim_bw_1}
\end{figure}

Figure~\ref{plot_sim_bw_1} shows the boxplots of our estimators $\hat{\eta}_{k_n}$ and $\hat{\ell}_{k_n,h_n}(\vec{x})$. In this example, the theoretical value of $\eta$ is $3/2+1=2.5$, and $\ell(\vec{x})$ is equal to $5.292757$ (cf. Table~\ref{tablecoeffs}). 

\subsection{Extreme risk measures estimation}

\begin{figure}[!h]
\begin{center}
\includegraphics[scale=0.6]{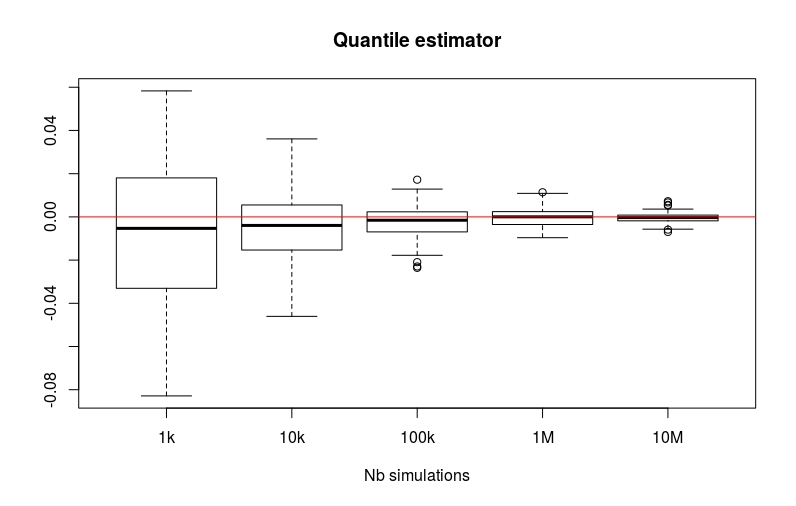}
\end{center}
\caption{Box plots representing $100$ relative errors $\frac{\hat{\hat{q}}_{\alpha_n} \left( Y| \vec{X}=\vec{x} \right)}{q_{\alpha_n} \left( Y| \vec{X}=\vec{x} \right)} - 1$ (based on sample sizes $n$ from 1 000 to 10 000 000) with $\alpha_n=1-n^{-1.25}$, $k_n=n^{0.6}$ and $h_n=n^{-0.2}$.}
\label{boxplots}
\end{figure}

It remains to estimate the conditional quantiles, expectiles and TVaRs of $Y| \vec{X}=\vec{x}$. Theoretical formulas (or algorithm) for conditional quantiles and expectiles may be found in~\citep{article} and~\cite{article2}. Furthermore, using straightforward calculations, formulas for Tail-Value-at-Risk may be obtained.
\begin{equation}
\left \{
\begin{array}{cc}
      q_{\alpha}\left(Y| \vec{X}=\vec{x} \right) =&  \sqrt{\frac{\nu+M(\vec{x})}{\nu+N}} \Phi_{\nu+N}^{-1}(\alpha)  \\
       \text{TVaR}_{\alpha}\left(Y| \vec{X}=\vec{x} \right) =& \frac{1}{1-\alpha} \frac{\Gamma \left( \frac{N+1+\nu}{2} \right)}{\Gamma \left( \frac{N+\nu}{2} \right)} \frac{\sqrt{\nu+M(\vec{x})}}{\sqrt{\pi} (\nu+N-1)} \left( 1+\frac{\Phi_{\nu+N}^{-1}(\alpha)^2}{\nu+N} \right)^{\frac{1-N-\nu}{2}} \\
\end{array},
\right.
\end{equation}
where $\Phi_{\nu}$ is the c.d.f of a Student distribution with $\nu$ degrees of freedom. In order to give an idea of the performance of our estimator, we propose in Figure~\ref{boxplots} some box plots representing $100$ relative errors (based on sample sizes $n$ from 1 000 to 10 000 000) of our quantile estimator~\eqref{eqqhathat} with $\alpha_n=1-n^{-1.25}$. \\ Finally, we would like to compare these results with other estimators already used. The most common and widespread methods for estimating conditional quantiles and expectiles are respectively quantile and expectile regression, introduced in~\cite{koenker} and~\cite{newey}. In~\cite{article} and~\cite{article2}, we have shown that such approach leads to a poor estimation in case of extreme levels. Indeed, in this example, a quantile regression estimator will converge to $\Phi_{\nu}^{-1}(\alpha_n)=1530.15$, very far from $7.31$, the theoretical result. Obviously, since the quantile regression estimator does not assume any structure on the underlying distribution, the latter is clearly less efficient than the tailored extreme quantile estimators introduced in this paper. 

It may also be interesting to compare the empirical variance of our estimator with our asymptotic result given in Proposition~\ref{propasympnormphirstar2}. Furthermore, the latter allow us to provide confidence intervals for $q_{\alpha_n}\left( Y| \vec{X}=\vec{x} \right)$. We thus introduce the notation $\hat{\zeta}_n$ the empirical variance of $\frac{\sqrt{k_n}}{\ln \left( \frac{k_n}{n (1-\alpha_n)} \right)} \left( \frac{\hat{\hat{q}}_{\alpha_n} \left(Y|\vec{X}=\vec{x} \right)}{q_{\alpha_n} \left(Y|\vec{X}=\vec{x} \right)}-1 \right)$, while $\zeta=\left( \frac{\gamma}{\gamma N+1}-\theta \frac{N \gamma^2}{(\gamma N+1)^2} \right)^2=0.0005536332$ in this section. In addition, we denote $m_n$ the number of times the theoretical value $q_{\alpha_n}\left( Y| \vec{X}=\vec{x} \right)$ is in the $95 \%$ confidence interval. Table~\ref{tableau_quantile} gives an overview of the behavior of these quantities according to $n$.

\begin{table}[!h]
\begin{tabular}{cccc}
  $n$ & $\alpha_n$ & $\hat{\zeta}_n$ & $m_n$ \\
  \hline
  1 000 & 0.9998222 & 0.001839654 & 44 \\
  10 000 & 0.99999 & 0.001110249 & 76 \\
  100 000 & 0.9999994 & 0.0006940025 & 90 \\
  1 000 000 & 0.99999996837 & 0.0005534661 & 94 \\
  10 000 000 & 0.99999999822 & 0.0005695426 & 91 \\
 $+\infty$ & 1 & 0.0005536332 & 95 \\
 \hline
\end{tabular}
\caption{Empirical variance $\hat{\zeta}_n$, number of confidence intervals containing the theoretical value $m_n$ for 100 estimations $\hat{\hat{q}}_{\alpha_n} \left( Y| \vec{X}=\vec{x} \right)$ of $q_{\alpha_n}\left( Y| \vec{X}=\vec{x} \right)$, with $n$ ranging from 1 000 to 10 000 000. Chosen sequences are $\alpha_n=1-n^{-1.25}$, $k_n=n^{0.6}$ and $h_n=n^{-0.2}$.}
\label{tableau_quantile}
\end{table}

\begin{figure}[!h]
\begin{center}
\includegraphics[scale=0.55]{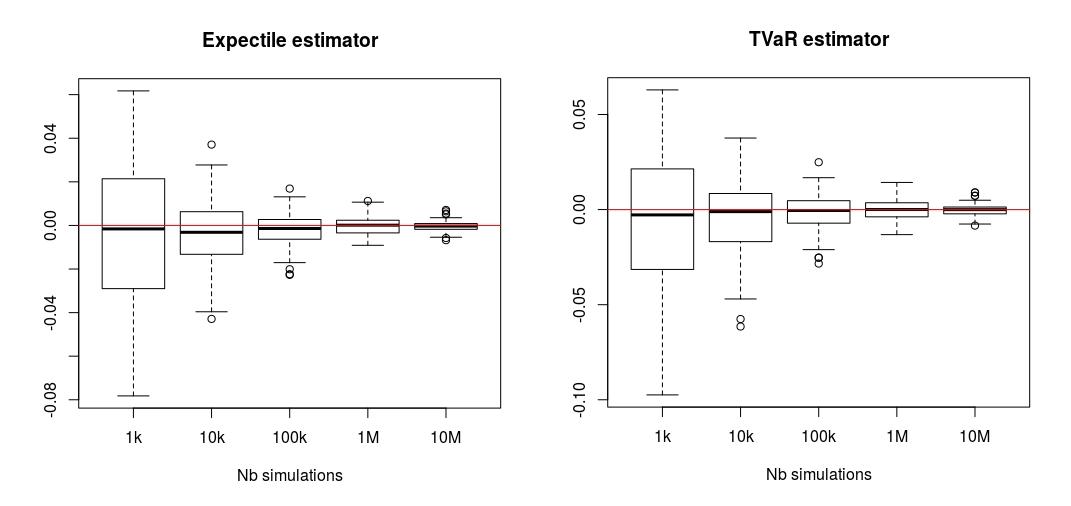}
\end{center}
\caption{From left to right : Box plots representing $100$ relative errors $\frac{\hat{\hat{q}}_{2,\alpha_n} \left( Y| \vec{X}=\vec{x} \right)}{q_{2,\alpha_n} \left( Y| \vec{X}=\vec{x} \right)} - 1$ and $\frac{\hat{\hat{H}}_{1,\alpha_n} \left( Y| \vec{X}=\vec{x} \right)}{H_{1,\alpha_n} \left( Y| \vec{X}=\vec{x} \right)} - 1$ (based on sample sizes $n$ from 1 000 to 10 000 000) with $\alpha_n=1-n^{-1.25}$, $k_n=n^{0.6}$ and $h_n=n^{-0.2}$.}
\label{boxplots_bis}
\end{figure}

Finally, based on these quantile estimates, we deduce, using Definitions~\ref{deflphat} and \ref{defhhat}, $L_2-$quantile (or expectile) and Tail-Value-at-Risk estimates. Figure~\ref{boxplots_bis} provides relative errors for estimators $\hat{\hat{q}}_{2,\alpha_n}$ and $\hat{\hat{H}}_{1,\alpha_n}$.

In the previous figures, only the first component of the vector is used to estimate the tail index. There is therefore some loss of information. We have suggested in Section~\ref{parametres} another approach. Furthermore, \cite{resnickdependant} or \cite{hsing} proved that the Hill estimator may also work with dependent data. Thus it would be possible to improve the estimation of $\hat{\gamma}_{k_n}$ by adding the other components of the vector in Equation~\eqref{hillest}, but in that case the asymptotic results of Propositions~\ref{propetahat} or \ref{proplhat} would not hold anymore. 

\section{Real data example} \label{real}

\begin{figure}[h!]
\begin{center}
\includegraphics[scale=0.4]{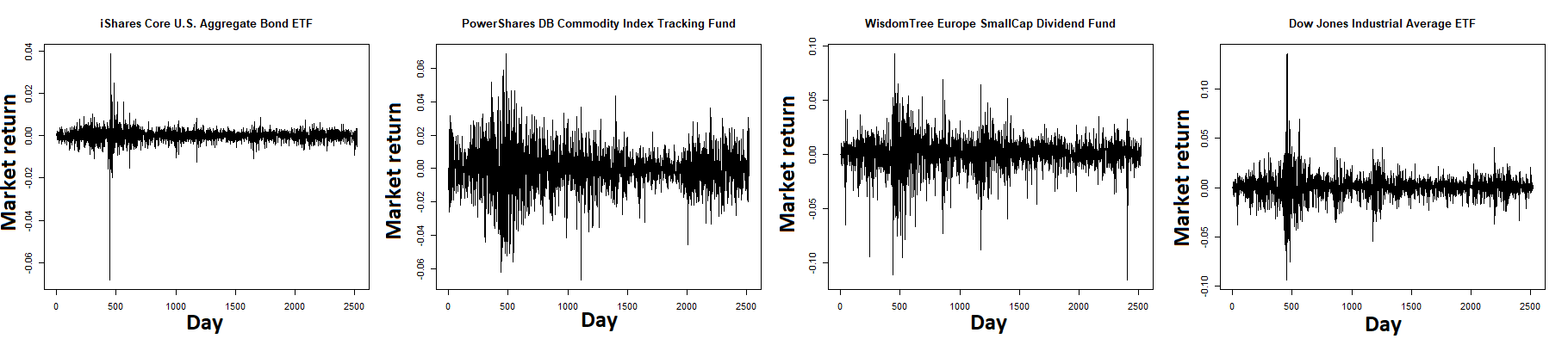}
\end{center}
\caption{Daily market returns of 4 different assets.}
\label{plotdata}
\end{figure}

As an application, we use the daily market returns (computed from the closing prices) of financial assets from 2006 to 2016, available at \url{http://stanford.edu/class/ee103/portfolio.html}. We focus on the first four assets, i.e iShares Core U.S. Aggregate Bond ETF, PowerShares DB Commodity Index Tracking Fund, WisdomTree Europe SmallCap Dividend Fund and SPDR Dow Jones Industrial Average ETF which will be our covariate $\vec{X}$. Figure~\ref{plotdata} represents the daily return for each day. 

The reason for focusing solely on the value of these assets could be, for example, that they are the first available every day. The aim would be to anticipate the behavior of another asset on another market. We thus consider the return of WisdomTree Japan Hedged Equity Fund as random variable $Y$. The size of the sample is 2520. The first 2519 days (from January 3, 2007 to December 5, 2016) will be our learning sample, and we focus on the 2520th day, when the covariate $\vec{X}$ is $\vec{x}=\left( -0.0185 \%,-0.4464 \%, 0.9614 \%, 0.1405 \% \right)$. Pending the opening of the second market, let us estimate the quantile of the return $Y$ given $\vec{X}=\vec{x}$. \\
After a brief study of the autocorrelation functions, we consider that the daily returns can be considered as independent. Concerning the shape of the data, histograms of the marginals seem symmetrical. Furthermore, the measured tail index is approximately the same for the 4 marginals. This is why suppose that the data is elliptical. After having estimated $\vec{\mu}$ and $\mat{\Sigma}$ by the method of moments, we get $M(\vec{x})=1.072952$. We apply our estimators $\hat{\eta}_{k_n}$ and $\hat{\ell}_{k_n,h_n}(\vec{x})$ given in Equations~\eqref{eqetahat} and~\eqref{eqlhat}. We take as sequences $k_n=n^{0.6}$ ($b=0.6$) and $h_n=n^{-0.2}$ ($c=0.2$), and as kernel $K$ the gaussian p.d.f, hence we deduce the asymptotic confidence bounds from Equation~\eqref{eqloijointe1}. We then obtain $\hat{\eta}_{k_n}=2.617846$ and $\hat{\ell}_{k_n,h_n}(\vec{x})=6.44334$. Let us now estimate the high quantile $q_{\alpha_n}\left( Y| \vec{X}=\vec{x} \right)$ with level $\alpha_n=1-n^{-a}$, $a>1$. In order to minimize the asymptotic variance of Equation~\eqref{eqasympnormphirstar2}, we chose $a=(1-b)\left( \hat{\gamma}_{k_n}+1 \right)=1.047146$. By applying estimator~\eqref{eqqhathat}, we get a quantile of level $0.9997256$ close to $3.744985 \%$ for $Y| \vec{X}=\vec{x}$. In other words, before the opening of the second market, we consider that given the returns of our first four assets, that of WisdomTree Japan Hedged Equity Fund has a probability $0.9997256$ of beeing less than $3.744985 \%$. For information, the true return that day was $0.7141 \%$.

\section{Conclusion} \label{Conclusion}

In this paper, we propose two estimators $\hat{\ell}_{k_n,h_n}(\vec{x})$ and $\hat{\eta}_{k_n}$ respectively for extremal parameters $\ell(\vec{x})$ and $\eta$ introduced in Equation~\eqref{ecqe}. We have proved their consistency and asymptotic normality according to the asymptotic relationships between the sequences $k_n$ and $h_n$. Using these estimators, we have defined estimators for intermediate and high quantiles, proved their consistency, given their asymptotic normality under stronger conditions, and deduced estimators for extreme $L_p-$quantiles and Haezendonck-Goovaerts risk measures. Consistency and asymptotic normality are also provided for these estimators, under conditions. We have also illustrated with a numerical example the performance of our estimators, and applied them to real data set. \\
As working perspectives, we intend to propose a method of optimal choice of the sequences $k_n$ and $h_n$, which is not totally discussed in this paper. Furthermore, the shape of $\ell(\vec{x})$ and $\eta$ leaving Assumption~\ref{hyp2} is a current research topic. More generally, the asymptotic relationships between conditional and unconditional quantile in other maximum domains of attraction, using for example the results of~\cite{hashorva}, may be developed. However, we need a second-order refinement, as we need a second-order refinement of Equation~\eqref{equivalencephirstar} to propose asymptotic normalities~\ref{propnormint2} and~\ref{propasympnormphirstar2} under weaker assumptions than Assumption~\ref{hypforte}. Finally, it seems that the ratio of the two terms in Equation~\eqref{equivalencephirstar} tends to 1 more and more slowly when the covariate vector size $N$ becomes large. Then, our estimation approach may perform poorly if $N$ is tall. This is why it might be wise to propose another method when the covariate vector size $N$ is large.

\section*{Appendix}

\subsection*{Proof of Lemma~\ref{lemmeprops}}

\begin{enumerate}[label=(\roman*)]
\item Since $R_1 \overset{d}{=} \chi_1 \xi$, where $\chi_1$ has a Lebesgue density $\sqrt{\frac{2}{\pi}} e^{-\frac{x^2}{2}}$. According to Lemma 4.3 in~\cite{jessen2006}, $\xi$ satisfies $\bar{F}_{\xi}(t \omega)/\bar{F}_{\xi}(t) \rightarrow \omega^{-\frac{1}{\gamma}} $ as $t \rightarrow +\infty$. Furthermore, Lemma 4.2 in~\cite{jessen2006} entails $$ \mathbb{P} \left( \xi >t \right) \underset{t \rightarrow +\infty}{\sim}\mathbb{E}\left[ \chi_1^{\frac{1}{\gamma}} \right]^{-1} \mathbb{P} \left( R_1>t \right).  $$ Assumption~\ref{hyp2} provides $\mathbb{P}(R_1>t) \sim \lambda t^{-\frac{1}{\gamma}}$, hence the result.
\item Using again Lemma 4.2 in~\cite{jessen2006} for $R_d \overset{d}{=} \chi_d \xi$, it comes immediatly $$ \mathbb{P}\left( R_d>t \right) \underset{t \rightarrow +\infty}{\sim} \mathbb{E}\left[ \chi_d^{\frac{1}{\gamma}} \right] \mathbb{P} \left( \xi >t \right).$$ Some straightforward calculations provide $\mathbb{E}\left[ \chi_d^{\frac{1}{\gamma}} \right]=2^{\frac{1}{\gamma}} \frac{\Gamma \left( \frac{d+\gamma^{-1}}{2} \right)}{\Gamma \left( \frac{d}{2} \right)} $.
\item From (ii), we have, for all $d \in \mathbb{N}$, $f_{R_d}(t)\underset{t \rightarrow +\infty}{\sim} 2^{\frac{1}{\gamma}} \frac{\Gamma \left( \frac{d+\gamma^{-1}}{2} \right)}{\Gamma \left( \frac{d}{2} \right)} \lambda' t^{-\frac{1}{\gamma}-1} $, where $\lambda' \in \mathbb{R}$ is not related to $d$. The result is immediate with this expression. $\Box$
\end{enumerate}  

\subsection*{Proof of Proposition~\ref{propetal}}

The conditional density (Proposition 3 in \cite{article}) leads to~:$$ \underset{t \rightarrow \infty}{\lim} \frac{\bar{\Phi}_{R^*}(t)}{\bar{\Phi}_{R}(t^{\eta})}= \underset{t \rightarrow \infty}{\lim} \frac{c_{N+1} g_{N+1}(M(\vec{x})+t^2)}{c_N g_N\left(M(\vec{x}) \right) \eta t^{\eta-1} c_1 g_1(t^{2 \eta})}=\underset{t \rightarrow \infty}{\lim} \frac{\Gamma \left( \frac{N+1}{2} \right)(M(\vec{x})+t^2)^{-\frac{N}{2}}}{\pi^{\frac{N+1}{2}} c_N g_N\left(M(\vec{x}) \right) \eta t^{\eta-1}} \frac{f_{R_{N+1}} \left( \sqrt{M(\vec{x})+t^2} \right)}{f_{R_{1}}(t^\eta)}.$$ Using Equation~\eqref{lemme1_3} of Lemma~\ref{lemmeprops}, it comes $$ \frac{\bar{\Phi}_{R^*}(t)}{\bar{\Phi}_{R}(t^{\eta})} \underset{t \rightarrow +\infty}{\sim}  \frac{1}{\pi^{\frac{N}{2}} c_N g_N\left(M(\vec{x}) \right) \eta} \frac{ \Gamma \left( \frac{N+1+\gamma^{-1}}{2} \right)}{\Gamma \left( \frac{1+\gamma^{-1}}{2} \right)} t^{(\eta-1)(\gamma^{-1}+1)+1-\eta-N}.  $$ Obviously, we impose $0<\ell(\vec{x})<+\infty$, then $1-\eta-N+(\eta-1)(\gamma^{-1}+1)=0$, hence $\eta=N\gamma+1$. Replacing $\eta$ in the previous equation, $\ell(\vec{x})$ is easily deduced~: $$ \ell(\vec{x})= \frac{1}{\pi^{\frac{N}{2}} c_N g_N\left(M(\vec{x}) \right) \eta} \frac{ \Gamma \left( \frac{N+1+\gamma^{-1}}{2} \right)}{\Gamma \left( \frac{1+\gamma^{-1}}{2} \right)}. \text{  } \Box $$

\subsection*{Proof of Proposition~\ref{propnorml}}

It is obvious that under conditions $(K1)-(K2)$, $\sqrt{k_n} \left( \hat{g}_{h_n}-c_N g_N\left(M(\vec{x}) \right) \right) \overset{\mathbb{P}}{\rightarrow} 0$ as $n \rightarrow +\infty$ if $k_n=o(n h_n)$ and $\sqrt{k_n}h_n^2 \rightarrow 0$. Then we get the following asymptotic normality~: $$ \sqrt{k_n} \begin{pmatrix}
   \hat{\gamma}_{k_n}-\gamma  \\
   \hat{g}_{h_n}-c_N g_N\left(M(\vec{x}) \right) 
\end{pmatrix} \underset{n \rightarrow +\infty}{\rightarrow} \mathcal{N} \left( \begin{pmatrix}
   0  \\
   0  
\end{pmatrix}, \begin{pmatrix}
   \gamma^2 & 0 \\
   0 & 0 
\end{pmatrix} \right).$$ Since $\ell(\vec{x})=u(\gamma)$, the delta method entails $$ \sqrt{k_n} \left( \hat{\ell}_{k_n,h_n}(\vec{x})-\ell(\vec{x}) \right) \underset{n \rightarrow +\infty}{\rightarrow} \mathcal{N} \left( 0,u'(\gamma)^2 \gamma^2 \right).$$ A quick calculation of $u'$, using Equation~\eqref{eqetal}, gives the first result. The second part of the proof is similar. Indeed, if $n h_n=o(k_n)$ and $n h_n^5 \rightarrow 0$ as $n \rightarrow +\infty$, then $$ \sqrt{n h_n} \begin{pmatrix}
   \hat{\gamma}_{k_n}-\gamma  \\
   \hat{g}_{h_n}-c_N g_N\left(M(\vec{x}) \right) 
\end{pmatrix} \underset{n \rightarrow +\infty}{\rightarrow} \mathcal{N} \left( \begin{pmatrix}
   0  \\
   0  
\end{pmatrix}, \begin{pmatrix}
   0 & 0 \\
   0 & \frac{M(\vec{x})^{1-\frac{N}{2}} \Gamma \left( \frac{N}{2} \right)}{\pi^{\frac{N}{2}}} c_N g_N \left(M(\vec{x}) \right) \int K(u)^2 du 
\end{pmatrix} \right).$$ The delta method completes the proof. $\Box$

In order to make the proof of Theorem~\ref{propnormint1} easier to read, we give the following lemma, which provides the asymptotic behavior of a statistic order under Assumption~\ref{hyp2}.

\begin{lemme}
Under Assumption~\ref{hyp2} and condition $(C)$,
\begin{equation}
\sqrt{n v_n} \left( \frac{ \left(W_{[n v_n+1]}\right)}{\Phi_R^{-1} \left( 1-v_n \right)}-1 \right) \underset{n \rightarrow +\infty}{\rightarrow}  \mathcal{N} \left( 0, \gamma^2  \right).
\end{equation}
\label{lemmeappendix}
\end{lemme} 

\subsection*{Proof of Lemma~\ref{lemmeappendix}}

The proof is inspired by Theorem 2.4.1 in \cite{dehaan2006}. Let $Y_1,Y_2,\hdots$ be independant and identically distributed random variables with c.d.f. $1-y^{-1}, y>1$. We denote in addition $Y_{[n]} \leq \hdots \leq Y_{[1]}$. We thus have $$ \sqrt{n v_n} \left( v_n Y_{[n v_n+1]}-1 \right) \underset{n \rightarrow +\infty}{\rightarrow} \mathcal{N} (0,1).$$ By noticing that $ W_{[nv_n+1]} \overset{d}{=} \Phi_R^{-1} \left( 1-1/Y_{[nv_n+1]} \right) $, it comes $$ \sqrt{n v_n} \left( \frac{ \left(W_{[n v_n+1]}\right)}{\Phi_R^{-1} \left( 1-v_n \right)}-1 \right) \overset{d}{=} \sqrt{n v_n} \left( \frac{ \Phi_R^{-1} \left( 1-1/Y_{[nv_n+1]} \right)}{\Phi_R^{-1} \left( 1-v_n \right)}-\left( v_n Y_{[nv_n+1]} \right)^{\gamma} \right)+\sqrt{n v_n} \left( \left( v_n Y_{[nv_n+1]} \right)^{\gamma}-1 \right).$$ The delta method entails that the second term tends to $\mathcal{N}(0,\gamma^2)$. Moreover, Assumption~\ref{hyp2} and $\sqrt{k_n} A \left( \frac{n}{k_n} \right) \rightarrow 0$ as $n \rightarrow +\infty$ ensure the asymptotic nullity of the first term. 

\subsection*{Proof of Theorem~\ref{propnormint1}}

In a first time, we can notice $\tilde{v}_n$ is related to $\hat{\ell}_{k_n,h_n}(\vec{x})$. Then, according to Proposition~\ref{propnorml}, (i) entails that we can deal with $v_n$ instead of $\tilde{v}_n$ in Equation~\eqref{eqnormint1}. Furthermore, we give the decomposition~:
\begin{multline} \nonumber
\frac{\sqrt{k_n}}{\ln \left( 1-\alpha_n \right)} \left( \frac{ \hat{q}_{\alpha_n} \left( Y|\vec{X}=\vec{x} \right)}{q_{\alpha_n \uparrow} \left( Y|\vec{X}=\vec{x} \right)}-1 \right) \underset{n \rightarrow +\infty}{\sim} \frac{\sqrt{k_n}}{\ln \left( 1-\alpha_n \right)} \left( \frac{ \left(W_{[n v_n+1]}\right)^{1/\hat{\eta}_{k_n}}}{\Phi_R^{-1} \left( 1-v_n \right)^{1/\eta}}-1 \right)= \\ \frac{\sqrt{k_n}}{\ln \left( 1-\alpha_n \right)} \left( \frac{ \left(W_{[n v_n+1]}\right)^{1/\hat{\eta}_{k_n}}}{\Phi_R^{-1} \left( 1-v_n \right)^{1/\hat{\eta}_{k_n}}}-1 \right) \Phi_R^{-1} \left( 1-v_n \right)^{1/\hat{\eta}_{k_n}-1/\eta} + \frac{\sqrt{k_n}}{\ln \left( 1-\alpha_n \right)} \left( \Phi_R^{-1} \left( 1-v_n \right)^{1/\hat{\eta}_{k_n}-1/\eta}-1 \right).
\end{multline}
Under Assumption~\ref{hyp2}, and according to Proposition~\ref{propetahat} and Theorem 2.4.1 in~\cite{dehaan2006} (with $(C)$), we have~: $$ \left \{
\begin{array}{cc}
      \sqrt{k_n} \left( \frac{1}{\hat{\eta}_{k_n}} - \frac{1}{\eta} \right) & \underset{n \rightarrow +\infty}{\rightarrow}  \mathcal{N} \left( 0, \frac{N^2 \gamma^2}{\left( \gamma N+1 \right)^4}  \right) \\
       \sqrt{n v_n} \left( \frac{ \left(W_{[n v_n+1]}\right)}{\Phi_R^{-1} \left( 1-v_n \right)}-1 \right) & \underset{n \rightarrow +\infty}{\rightarrow}  \mathcal{N} \left( 0, \gamma^2  \right) \\
\end{array}
\right. , $$ By noticing that $v_n$ is equivalent to $\ell(\vec{x})^{-1}(1-\alpha_n)$ as $n \rightarrow +\infty$, and using condition $\sqrt{k_n}=o \left( \ln(1-\alpha_n) \sqrt{n(1-\alpha_n)} \right)$ in $(C_{int})$, it comes $$\frac{\sqrt{k_n}}{\ln \left( 1-\alpha_n \right)} \left( \frac{ \left(W_{[n v_n+1]}\right)^{1/\hat{\eta}_{k_n}}}{\Phi_R^{-1} \left( 1-v_n \right)^{1/\hat{\eta}_{k_n}}}-1 \right) \rightarrow 0 \text{ as } n \rightarrow +\infty. $$ Furthermore, under Assumption~\ref{hyp2}, $\ln \left( \Phi_R^{-1}(1-v_n) \right)$ is cleary equivalent to $-\gamma \ln(v_n)$, or $-\gamma \ln(1-\alpha_n)$. Then $(C_{int})$ ensures $\Phi_R^{-1} \left( 1-v_n \right)^{1/\hat{\eta}_{k_n}-1/\eta} \rightarrow 1$ as $n \rightarrow +\infty$, and therefore the first term of the decomposition tends to 0. It thus remains to calculate the limit of the second term. It is not complicated to notice that $$  \frac{\sqrt{k_n}}{\ln \left( \Phi_R^{-1} \left( 1-v_n \right) \right)} \left( \Phi_R^{-1} \left( 1-v_n \right)^{1/\hat{\eta}_{k_n}-1/\eta}-1 \right) \underset{n \rightarrow +\infty}{\rightarrow} \mathcal{N} \left( 0, \frac{N^2 \gamma^2}{\left( \gamma N+1 \right)^4}  \right).$$ Using the equivalence $\ln \left( \Phi_R^{-1}(1-v_n) \right) \sim -\gamma \ln(v_n) \sim -\gamma \ln(1-\alpha_n)$, we get the result~\eqref{eqnormint1}. Using asymptotic relationship~\eqref{equivalencephirstar}, the consistency~\ref{eqcvprobainterm} is obvious. $\Box$

\subsection*{Proof of Proposition~\ref{propnormint2}}

We recall that density of $\Phi_{R^*}$ is proportional to $c_{N+1} g_{N+1}\left( M(\vec{x})+t^2 \right)$, and, from Assumption~\ref{hypforte}, there exist $\lambda_1,\lambda_2 \in \mathbb{R}$ such that~:$$ c_{N+1} g_{N+1}\left( M(\vec{x})+t^2 \right)=\lambda_1 \left( M(\vec{x})+t^2 \right)^{-\frac{N+1+\gamma^{-1}}{2}} \left[ 1+\lambda_2 \left( M(\vec{x})+t^2 \right)^{\frac{\rho}{2\gamma}}+o \left( t^{\frac{\rho}{\gamma}}  \right) \right] $$ The previous expression may be rewritten as follows, where $\lambda_1,\lambda_2,\lambda_3 \in \mathbb{R}$~: $$ c_{N+1} g_{N+1}\left( M(\vec{x})+t^2 \right)=\lambda_1 t^{-(N+1+\gamma^{-1})} \left[ 1+\lambda_2 \left( M(\vec{x})+t^2 \right)^{\frac{\rho}{2\gamma}}+\lambda_3 t^{-2} +o \left( t^{\frac{\rho}{\gamma}}  \right) \right] $$ In order to make the proof more readable, we do not specify the values of constants $\lambda_i$, because they are not essential. Then, in the case, $\rho/\gamma \leq -2$, we get $$ c_{N+1} g_{N+1}\left( M(\vec{x})+t^2 \right)=\lambda_1 t^{-(N+1+\gamma^{-1})} \left[ 1+\lambda_2 t^{-2} +o \left( t^{-2}  \right) \right], \lambda_1,\lambda_2 \in \mathbb{R} $$ In other terms, $c_{N+1} g_{N+1}\left( M(\vec{x})+t^2 \right)$ is regularly varying of second order with indices $-N-1-\gamma^{-1}$, $-2$, and an auxiliary function propotional to $t^{-2}$. According to Proposition 6 of~\cite{hua}, $\bar{\Phi}_{R^*}(t)=\int_t^{+\infty} c_{N+1} g_{N+1}\left( M(\vec{x})+u^2 \right) du \in 2RV_{-N-\gamma^{-1},-2}$ with an auxiliary function proportional to $t^{-2}$. Equivalently, there exists $\lambda_1,\lambda_2 \in \mathbb{R}$ such that $$ \Phi_{R^*}^{-1} \left( 1-\frac{1}{t} \right)=\lambda_1 t^{\frac{\gamma}{\gamma N+1}} \left[ 1+\lambda_2 t^{-\frac{2 \gamma}{\gamma N+1}}+o \left( t^{-\frac{2 \gamma}{\gamma N+1}} \right) \right]. $$ Since Assumption~\ref{hyp2} and Assumption~\ref{hypforte} provide $\Phi_R^{-1}(1-1/t)=\lambda_3 t^{\gamma} \left[ 1+\lambda_4 t^{\rho}+o \left( t^{\rho} \right) \right]$, it comes $$ \frac{\Phi_{R}^{-1} \left( 1-v_n \right)^{\frac{1}{\eta}}}{\Phi_{R^*}^{-1}(\alpha_n)}=\ell(\vec{x})^{-\frac{\gamma}{\gamma N+1}} \left( \frac{1-\alpha_n}{v_n} \right)^{\frac{\gamma}{\gamma N+1}} \frac{1+\lambda_1 v_n^{-\rho}+o\left( v_n^{-\rho} \right)}{1+\lambda_2 (1-\alpha_n)^{\frac{2 \gamma}{\gamma N+1}} +o\left( (1-\alpha_n)^{\frac{2 \gamma}{\gamma N+1}} \right)},  $$ for some constants $\lambda_1,\lambda_2 \in \mathbb{R}$. In that case, we considered $\rho \leq -2\gamma$, hence $-\rho>2\gamma/(\gamma N+1)$. We then deduce the following expansion~: $$ \frac{\Phi_{R}^{-1} \left( 1-v_n \right)^{\frac{1}{\eta}}}{\Phi_{R^*}^{-1}(\alpha_n)}=\ell(\vec{x})^{-\frac{\gamma}{\gamma N+1}} \left( \frac{1-\alpha_n}{v_n} \right)^{\frac{\gamma}{\gamma N+1}} \left[ 1+\lambda (1-\alpha_n)^{\frac{2 \gamma}{\gamma N+1}} +o\left( (1-\alpha_n)^{\frac{2 \gamma}{\gamma N+1}} \right) \right], $$ for a certain constant $\lambda \in \mathbb{R}$. We can notice that $(1-\alpha_n)/v_n=2(1-\ell(\vec{x}))(1-\alpha_n)+\ell(\vec{x})$, and let us now focus on the limit~:
\begin{multline}\nonumber
\underset{n \rightarrow +\infty}{\lim} \text{ } \frac{\sqrt{k_n}}{\ln \left( 1-\alpha_n \right)} \ln \left( \frac{\Phi_{R}^{-1} \left( 1-v_n \right)^{\frac{1}{\eta}}}{\Phi_{R^*}^{-1}(\alpha_n)} \right)= \frac{\gamma}{\gamma N+1} \underset{n \rightarrow +\infty}{\lim} \text{ } \frac{\sqrt{k_n}}{\ln \left( 1-\alpha_n \right)} \ln \left( 2 \frac{1-\ell(\vec{x})}{\ell(\vec{x})} (1-\alpha_n)+1 \right) \\ + \underset{n \rightarrow +\infty}{\lim} \text{ } \frac{\sqrt{k_n}}{\ln \left( 1-\alpha_n \right)} \ln \left( 1+\lambda (1-\alpha_n)^{\frac{2 \gamma}{\gamma N+1}} +o\left( (1-\alpha_n)^{\frac{2 \gamma}{\gamma N+1}} \right) \right).
\end{multline}
The first term gives is easy to calculate. Indeed, since $\sqrt{k_n}(1-\alpha_n)/\ln \left( 1-\alpha_n \right) \rightarrow 0$ as $n \rightarrow +\infty$, we deduce $$ \underset{n \rightarrow +\infty}{\lim} \text{ } \frac{\sqrt{k_n}}{\ln \left( 1-\alpha_n \right)} \ln \left( 2 \frac{1-\ell(\vec{x})}{\ell(\vec{x})} (1-\alpha_n)+1 \right)=2 \frac{1-\ell(\vec{x})}{\ell(\vec{x})} \underset{n \rightarrow +\infty}{\lim} \text{ } \frac{\sqrt{k_n}}{\ln \left( 1-\alpha_n \right)} (1-\alpha_n)=0.$$ By a similar calculation, the second term also tends to 0, supposing $\frac{\sqrt{k_n}}{\ln \left( 1-\alpha_n \right)} \left( 1-\alpha_n \right)^{\frac{2 \gamma}{\gamma N+1}} \rightarrow 0$ as $n \rightarrow +\infty$. Then, we deduce, using Proposition~\ref{propnormint1}~:
\begin{multline} \nonumber
\frac{\sqrt{k_n}}{\ln \left( 1-\alpha_n \right)} \left( \frac{\hat{q}_{\alpha_n} \left( Y| \vec{X}=\vec{x} \right)}{q_{\alpha_n} \left( Y| \vec{X}=\vec{x} \right)}-1 \right) \sim \frac{\sqrt{k_n}}{\ln \left( 1-\alpha_n \right)} \left( \frac{ \left(W_{[n v_n+1]}\right)^{1/\hat{\eta}_{k_n}}}{\Phi_{R^*}^{-1} \left( \alpha_n \right)}-1 \right)= \\ \frac{\sqrt{k_n}}{\ln \left( 1-\alpha_n \right)} \left( \frac{\hat{q}_{\alpha_n} \left(R^* U^{(1)} \right)}{\Phi_{R}^{-1} \left( 1-v_n \right)^{\frac{1}{\eta}}}-1 \right) \frac{\Phi_{R}^{-1} \left( 1-v_n \right)^{\frac{1}{\eta}}}{\Phi_{R^*}^{-1}(\alpha_n)} + \frac{\sqrt{k_n}}{\ln \left( 1-\alpha_n \right)} \left( \frac{\Phi_{R}^{-1} \left( 1-v_n \right)^{\frac{1}{\eta}}}{\Phi_{R^*}^{-1}(\alpha_n)}-1 \right) \underset{n \rightarrow +\infty}{\rightarrow} \mathcal{N} \left( 0,\frac{N^2 \gamma^4}{\left( \gamma N+1 \right)^4} \right).
\end{multline}
Now, let us focus on the case $\rho/\gamma > -2$. The proof is exactly the same, with $$ c_{N+1} g_{N+1}\left( M(\vec{x})+t^2 \right)=\lambda_1 t^{-(N+1+\gamma^{-1})} \left[ 1+\lambda_2 t^{\frac{\rho}{\gamma}} +o \left( t^{\frac{\rho}{\gamma}}  \right) \right], \lambda_1,\lambda_2 \in \mathbb{R}. $$ Using the same calculations and doing the further assumption $\underset{n \rightarrow +\infty}{\lim} \frac{\sqrt{k_n}}{\ln \left( 1-\alpha_n \right)} \left( 1-\alpha_n \right)^{-\frac{\rho}{\gamma N+1}}=0$ leads to the result. $\Box$

\subsection*{Proof of Theorem~\ref{thasympnormphirstar}}

Firstly, we can notice that 
\begin{multline}
\frac{\hat{\hat{q}}_{\alpha_n} \left( Y| \vec{X}=\vec{x} \right)}{q_{\alpha_n \uparrow} \left( Y| \vec{X}=\vec{x} \right)}-1 \underset{n \rightarrow +\infty}{\sim} \frac{\left[ W_{[k_n+1]} \left( \frac{k_n}{n \tilde{v}_n} \right)^{\hat{\gamma}_{k_n}} \right]^{\frac{1}{\hat{\eta}_{k_n}}}}{\Phi_R^{-1} \left( 1-v_n \right)^{\frac{1}{\eta}}}-1 = \\ \left( \frac{\left[ W_{[k_n+1]} \left( \frac{k_n}{n v_n} \right)^{\hat{\gamma}_{k_n}} \right]^{\frac{1}{\hat{\eta}_{k_n}}}}{\Phi_R^{-1} \left( 1-v_n \right)^{\frac{1}{\eta}}}-1 \right) \left( \frac{v_n}{\tilde{v}_n} \right)^{\frac{\hat{\gamma}_{k_n}}{\hat{\eta}_{k_n}}}+\left( \frac{v_n}{\tilde{v}_n} \right)^{\frac{\hat{\gamma}_{k_n}}{\hat{\eta}_{k_n}}}-1.
\end{multline}
Since $k_n=o(nh_n)$, we deduce $$\frac{\sqrt{k_n}}{\ln \left( \frac{k_n}{n v_n} \right)} \left(  \frac{\left[ W_{[k_n+1]} \left( \frac{k_n}{n \tilde{v}_n} \right)^{\hat{\gamma}_{k_n}} \right]^{\frac{1}{\hat{\eta}_{k_n}}}}{\Phi_R^{-1} \left( 1-v_n \right)^{\frac{1}{\eta}}}-1 \right) \underset{n \rightarrow +\infty}{\sim} \frac{\sqrt{k_n}}{\ln \left( \frac{k_n}{n v_n} \right)} \left(  \frac{\left[ W_{[k_n+1]} \left( \frac{k_n}{n v_n} \right)^{\hat{\gamma}_{k_n}} \right]^{\frac{1}{\hat{\eta}_{k_n}}}}{\Phi_R^{-1} \left( 1-v_n \right)^{\frac{1}{\eta}}}-1 \right). $$ Furthermore, according to Theorem 4.3.8 in~\cite{dehaan2006}, $(C)$ and $(C_{high})$ lead to $$ \frac{\sqrt{k_n}}{\ln \left( \frac{k_n}{n v_n} \right)} \left(  \frac{ W_{[k_n+1]} \left( \frac{k_n}{n v_n} \right)^{\hat{\gamma}_{k_n}} }{\Phi_R^{-1} \left( 1-v_n \right)}-1 \right)\underset{n \rightarrow +\infty}{\sim} \frac{\sqrt{k_n}}{\ln \left( \frac{k_n}{n v_n} \right)} \left( \left( \frac{k_n}{n v_n} \right)^{\hat{\gamma}_{k_n}-\gamma}-1  \right). $$ From Assumption~\ref{hyp2}, it is not difficult to prove that $\ln\left( \Phi_R^{-1}(1-v_n) \right)/\ln \left( k_n/(n v_n) \right) $ is asymptotically equivalent to $ \gamma \ln(1-\alpha_n)/ \ln \left( n(1-\alpha_n)/k_n \right)$. Then, if we focus on the second term, it comes, using the limit given in $(C_{high})$~: $$ \frac{\sqrt{k_n}}{\ln \left( \frac{k_n}{n v_n} \right)} \begin{pmatrix}
    \left( \frac{k_n}{n v_n} \right)^{\hat{\gamma}_{k_n}-\gamma}-1 \\
   \Phi_R\left( 1-v_n \right)^{ \frac{1}{\hat{\eta}_{k_n}} -\frac{1}{\eta}}-1 
\end{pmatrix} \underset{n \rightarrow +\infty}{\rightarrow} \mathcal{N} \left( \begin{pmatrix}
   0  \\
   0  
\end{pmatrix}, \begin{pmatrix}
   \gamma^2 & -\theta \frac{N \gamma^3}{(\gamma N+1)^2} \\
   -\theta \frac{N \gamma^3}{(\gamma N+1)^2} & \theta^2 \frac{N^2 \gamma^4}{(\gamma N+1)^4} 
\end{pmatrix} \right).  $$ Finally,
\begin{multline} 
\frac{\sqrt{k_n}}{\ln \left( \frac{k_n}{n v_n} \right)} \left( \frac{\left[ W_{[k_n+1]} \left( \frac{k_n}{n \tilde{v}_n} \right)^{\hat{\gamma}_{k_n}} \right]^{\frac{1}{\hat{\eta}_{k_n}}}}{\Phi_R^{-1} \left( 1-v_n \right)^{\frac{1}{\eta}}}-1 \right)=  \frac{\sqrt{k_n}}{\ln \left( \frac{k_n}{n v_n} \right)}  \left( \Phi_R^{-1} \left( 1-v_n \right)^{\frac{1}{\hat{\eta}_{k_n}}-\frac{1}{\eta}}-1  \right) \\ + \frac{\sqrt{k_n}}{\ln \left( \frac{k_n}{n v_n} \right)} \left( \frac{\left[ W_{[k_n+1]} \left( \frac{k_n}{n \tilde{v}_n} \right)^{\hat{\gamma}_{k_n}} \right]^{\frac{1}{\hat{\eta}_{k_n}}}}{\Phi_R^{-1} \left( 1-v_n \right)^{\frac{1}{\hat{\eta}_{k_n}}}}-1 \right) \Phi_R^{-1} \left( 1-v_n \right)^{\frac{1}{\hat{\eta}_{k_n}}-\frac{1}{\eta}}.
\end{multline}
When $n \rightarrow \infty$, this expression is the sum of the following bivariate normal distribution~: $$ \frac{\sqrt{k_n}}{\ln \left( \frac{k_n}{n v_n} \right)} \begin{pmatrix}
    \frac{\left[ W_{[k_n+1]} \left( \frac{k_n}{n \tilde{v}_n} \right)^{\hat{\gamma}_{k_n}} \right]^{\frac{1}{\hat{\eta}_{k_n}}}}{\Phi_R^{-1} \left( 1-v_n \right)^{\frac{1}{\hat{\eta}_{k_n}}}}-1  \\
   \Phi_R\left( 1-v_n \right)^{ \frac{1}{\hat{\eta}_{k_n}} -\frac{1}{\eta}}-1 
\end{pmatrix} \underset{n \rightarrow +\infty}{\rightarrow} \mathcal{N} \left( \begin{pmatrix}
   0  \\
   0  
\end{pmatrix}, \begin{pmatrix}
    \frac{\gamma^2}{(\gamma N+1)^2} & -\theta \frac{N \gamma^3}{(\gamma N+1)^3} \\
   -\theta \frac{N \gamma^3}{(\gamma N+1)^3} & \theta^2 \frac{N^2 \gamma^4}{(\gamma N+1)^4} 
\end{pmatrix} \right),  $$ To conclude, $\ln \left( \frac{k_n}{n v_n} \right) \sim \ln \left( \frac{k_n}{n (1-\alpha_n)} \right)$ as $n \rightarrow +\infty$, hence the result. The consistency is immediate. $\Box$

\subsection*{Proof of Proposition~\ref{propasympnormphirstar2}}

The proof is similar to that of Proposition~\ref{propnormint2}. Indeed, we have given, in the case $\rho/\gamma \leq -2$~: $$ \frac{\Phi_{R}^{-1} \left( 1-v_n \right)^{\frac{1}{\eta}}}{\Phi_{R^*}^{-1}(\alpha_n)}=\ell(\vec{x})^{-\frac{\gamma}{\gamma N+1}} \left( 2(1-\ell(\vec{x}))(1-\alpha_n)+\ell(\vec{x}) \right)^{\frac{\gamma}{\gamma N+1}} \left[ 1+\lambda (1-\alpha_n)^{\frac{2 \gamma}{\gamma N+1}} +o\left( (1-\alpha_n)^{\frac{2 \gamma}{\gamma N+1}} \right) \right], $$ for a certain constant $\lambda \in \mathbb{R}$. It thus remains to calculate
\begin{multline}\nonumber
\underset{n \rightarrow +\infty}{\lim} \text{ } \frac{\sqrt{k_n}}{\ln \left( \frac{k_n}{n (1-\alpha_n)} \right)} \ln \left( \frac{\Phi_{R}^{-1} \left( 1-v_n \right)^{\frac{1}{\eta}}}{\Phi_{R^*}^{-1}(\alpha_n)} \right)= \frac{\gamma}{\gamma N+1} \underset{n \rightarrow +\infty}{\lim} \text{ } \frac{\sqrt{k_n}}{\ln \left( \frac{k_n}{n (1-\alpha_n)} \right)} \ln \left( 2 \frac{1-\ell(\vec{x})}{\ell(\vec{x})} (1-\alpha_n)+1 \right) \\ + \underset{n \rightarrow +\infty}{\lim} \text{ } \frac{\sqrt{k_n}}{\ln \left( \frac{k_n}{n (1-\alpha_n)} \right)} \ln \left( 1+\lambda (1-\alpha_n)^{\frac{2 \gamma}{\gamma N+1}} +o\left( (1-\alpha_n)^{\frac{2 \gamma}{\gamma N+1}} \right) \right).
\end{multline}
The first term is easy to calculate. Indeed, since $n(1-\alpha_n) \rightarrow 0$ and $k_n=o(n)$ as $n \rightarrow +\infty$, we deduce $$ \underset{n \rightarrow +\infty}{\lim} \text{ } \frac{\sqrt{k_n}}{\ln \left( \frac{k_n}{n (1-\alpha_n)} \right)} \ln \left( 2 \frac{1-\ell(\vec{x})}{\ell(\vec{x})} (1-\alpha_n)+1 \right)=2 \frac{1-\ell(\vec{x})}{\ell(\vec{x})} \underset{n \rightarrow +\infty}{\lim} \text{ } \frac{\sqrt{k_n}}{\ln \left( \frac{k_n}{n (1-\alpha_n)} \right)} (1-\alpha_n)=0.$$ By a similar calculation, the second term also tends to 0, supposing $\frac{\sqrt{k_n}}{\ln \left( \frac{k_n}{n (1-\alpha_n)} \right)} \left( 1-\alpha_n \right)^{\frac{2 \gamma}{\gamma N+1}} \rightarrow 0$ as $n \rightarrow +\infty$. Then, we deduce, using Proposition~\ref{thasympnormphirstar}~:
\begin{multline}\nonumber
\frac{\sqrt{k_n}}{\ln \left( \frac{k_n}{n (1-\alpha_n)} \right)} \left( \frac{\left[ W_{[k_n+1]} \left( \frac{k_n}{n \tilde{v}_n} \right)^{\hat{\gamma}_{k_n}} \right]^{\frac{1}{\hat{\eta}_{k_n}}}}{\Phi_{R^*}^{-1} \left( \alpha_n \right)}-1 \right)=\frac{\sqrt{k_n}}{\ln \left( \frac{k_n}{n (1-\alpha_n)} \right)} \left( \frac{\left[ W_{[k_n+1]} \left( \frac{k_n}{n \tilde{v}_n} \right)^{\hat{\gamma}_{k_n}} \right]^{\frac{1}{\hat{\eta}_{k_n}}}}{\Phi_{R}^{-1} \left( 1-v_n \right)^{\frac{1}{\eta}}}-1 \right) \frac{\Phi_{R}^{-1} \left( 1-v_n \right)^{\frac{1}{\eta}}}{\Phi_{R^*}^{-1}(\alpha_n)} \\ + \frac{\sqrt{k_n}}{\ln \left( \frac{k_n}{n (1-\alpha_n)} \right)} \left( \frac{\Phi_{R}^{-1} \left( 1-v_n \right)^{\frac{1}{\eta}}}{\Phi_{R^*}^{-1}(\alpha_n)}-1 \right) \underset{n \rightarrow +\infty}{\rightarrow} \mathcal{N} \left( 0,\frac{\gamma^2}{(\gamma N+1)^2}-2\theta \frac{N \gamma^3}{(\gamma N+1)^3}+\theta^2 \frac{N^2 \gamma^4}{(\gamma N+1)^4} \right).
\end{multline}
Now, let us focus on the case $\rho/\gamma > -2$. The proof is exactly the same, with $$\frac{\Phi_{R}^{-1} \left( 1-v_n \right)^{\frac{1}{\eta}}}{\Phi_{R^*}^{-1}(\alpha_n)}=\ell(\vec{x})^{-\frac{\gamma}{\gamma N+1}} \left( 2(1-\ell(\vec{x}))(1-\alpha_n)+\ell(\vec{x}) \right)^{\frac{\gamma}{\gamma N+1}} \left[ 1+\lambda (1-\alpha_n)^{\frac{-\rho}{\gamma N+1}} +o\left( (1-\alpha_n)^{\frac{-\rho}{\gamma N+1}} \right) \right], \lambda \in \mathbb{R}. $$ Using the same calculations and doing the further assumption $\underset{n \rightarrow +\infty}{\lim} \frac{\sqrt{k_n}}{\ln \left( \frac{k_n}{n (1-\alpha_n)} \right)} \left( 1-\alpha_n \right)^{-\frac{\rho}{\gamma N+1}}=0$ leads to the result. $\Box$

\subsection*{Proof of Lemma~\ref{lemmegammacond}}

The density of $Y|\vec{X}=\vec{x}$ is given by $$ c_{N+1} g_{N+1} \left( M(\vec{x})+ (t-\mu_{Y|\vec{X}})^2 \sigma_{Y|\vec{X}}^{-2} \right) \left(c_N g_N \left( M(\vec{x}) \right) \right)^{-1},  $$ where $M(\vec{x})=(\vec{x}-\vec{\mu}_{\vec{X}})^{\top} \mat{\Sigma}_{\vec{X}}^{-1}(\vec{x}-\vec{\mu}_{\vec{X}})$. In order to simplify, we consider the case reduced and centered, i.e $\mu_{Y|\vec{X}}=0$ and $\sigma_{Y|\vec{X}}=1$. A quick calculation gives
$$ \underset{t \rightarrow +\infty}{\lim} \text{ } \frac{\bar{\Phi}_{R^*}(\omega t)}{\bar{\Phi}_{R^*}(t)}=\omega \underset{t \rightarrow +\infty}{\lim} \text{ } \frac{g_{N+1}(M(\vec{x})+\omega^2 t^2)}{g_{N+1}(M(\vec{x})+ t^2)} =\omega \underset{t \rightarrow +\infty}{\lim} \text{ } \frac{(M(\vec{x})+\omega^2 t^2)^{-\frac{N}{2}}}{(M(\vec{x})+ t^2)^{-\frac{N}{2}}} \frac{f_{R_{N+1}} \left( \sqrt{M(\vec{x})+ \omega^2 t^2} \right)}{f_{R_{N+1}} \left( \sqrt{M(\vec{x})+ t^2} \right)}.$$ Equation~\eqref{lemme1_1} leads to $$ \underset{t \rightarrow +\infty}{\lim} \frac{\bar{\Phi}_{R^*}(\omega t)}{\bar{\Phi}_{R^*}(t)}=\omega \omega^{-N} \omega^{-\frac{1}{\gamma}-1}=\omega^{-\frac{1}{\gamma}-N}. \text{   } \Box $$

\subsection*{Proof of Proposition~\ref{propnormlp}}

We recall in a first time that condition $\left( C_{int}^{L_p} \right)$ entails $\left( C_{int}^{HG} \right)$. We have the following decomposition~: 
\begin{multline}\nonumber
\frac{\sqrt{k_n}}{\ln \left( 1-\alpha_n \right)} \left( \frac{\hat{q}_{p,\alpha_n}\left(Y| \vec{X}=\vec{x} \right)}{q_{p,\alpha_n}\left(Y| \vec{X}=\vec{x} \right)}-1 \right)  \underset{n \rightarrow +\infty}{\sim} \frac{\sqrt{k_n}}{\ln \left( 1-\alpha_n \right)} \left( \frac{\left( W_{\left[n \tilde{v}_n +1 \right]} \right)^{1/\hat{\eta}_{k_n}} f_L \left( \left( \hat{\gamma}_{k_n}^{-1}+N \right)^{-1}, p \right)}{q_{p,\alpha_n}\left(R^* U^{(1)} \right)}-1 \right) = \\ \frac{\sqrt{k_n}}{\ln \left( 1-\alpha_n \right)} \left( \frac{f_L \left( \left( \hat{\gamma}_{k_n}^{-1}+N \right)^{-1}, p \right)}{f_L \left( \left( \gamma^{-1}+N \right)^{-1}, p \right)}-1  \right) \frac{\left( W_{\left[n \tilde{v}_n +1 \right]} \right)^{1/\hat{\eta}_{k_n}}}{\Phi_{R^*}^{-1}(\alpha_n)} \frac{f_L \left( \left( \gamma^{-1}+N \right)^{-1}, p \right) \Phi_{R^*}^{-1}(\alpha_n)}{q_{p,\alpha_n}\left(R^* U^{(1)} \right)} + \\ \frac{\sqrt{k_n}}{\ln \left( 1-\alpha_n \right)} \left( \frac{\left( W_{\left[n \tilde{v}_n +1 \right]} \right)^{1/\hat{\eta}_{k_n}}}{\Phi_{R^*}^{-1}(\alpha_n)}-1 \right) \frac{f_L \left( \left( \gamma^{-1}+N \right)^{-1}, p \right) \Phi_{R^*}^{-1}(\alpha_n)}{q_{p,\alpha_n}\left(R^* U^{(1)} \right)}+ \\ \frac{\sqrt{k_n}}{\ln \left( 1-\alpha_n \right)} \left( \frac{f_L \left( \left( \gamma^{-1}+N \right)^{-1}, p \right) \Phi_{R^*}^{-1}(\alpha_n)}{q_{p,\alpha_n}\left(R^* U^{(1)} \right)}-1 \right).
\end{multline}
We know that $f_L \left( \left( \hat{\gamma}_{k_n}^{-1}+N \right)^{-1}, p \right)$, as a function of $\hat{\gamma}_{k_n}$, is asymptotically normal with rate $\sqrt{k_n}$ (see Equation~\eqref{eqnormgammahat}). Then, the first term in the sum clearly tends to 0 as $n \rightarrow +\infty$. Using Proposition~\ref{propnormint2}, the second term tends to the normal distribution given in~\eqref{eqnormint2}. Finally, we have to check that the third term tends to 0. For that purpose, we use the second order expansion given in~\cite{girardlp}: $$ \frac{q_{p,\alpha_n}\left(R^* U^{(1)} \right)}{f_L \left( \left( \gamma^{-1}+N \right)^{-1}, p \right) q_{\alpha_n}\left( R^* U^{(1)} \right)}=1- (\gamma^{-1}+N)^{-1} r(\alpha_n,p)+ \left( \lambda+o(1) \right) A^* \left( \frac{1}{1-\alpha_n} \right),$$ where $r(\alpha_n,p)=\lambda_1 \frac{1}{q_{\alpha_n}\left( R^* U^{(1)} \right)} \left( \mathbb{E} \left[ R^* U^{(1)}  \right] +o(1) \right)+\lambda_2 A^* \left( \frac{1}{1-\alpha_n} \right) (1+o(1))$, $\lambda,\lambda_1,\lambda_2 \in \mathbb{R}$ are not related to $n$ and $A^*(t)$ is the auxiliary function of $\Phi_{R^*}\left( 1-\frac{1}{t} \right)$. It seems important to precise that the conditional distribution $R^* U^{(1)}$ is regularly varying with tail index $\gamma^{-1}+N > 1$, then $\mathbb{E} \left[ R^* U^{(1)}  \right]$ exists and, $R^* U^{(1)}$ being symmetric, equals 0. Then, a sufficient condition for asymptotic normality may be $$ \left \{
\begin{array}{cc}
      \underset{n \rightarrow +\infty}{\lim} \text{ } \frac{\sqrt{k_n}}{\ln \left( 1-\alpha_n \right) q_{\alpha_n}\left( R^* U^{(1)} \right)}  =& 0 \\
       \underset{n \rightarrow +\infty}{\lim} \text{ } \frac{\sqrt{k_n}}{\ln \left( 1-\alpha_n \right)} A^* \left( \frac{1}{1-\alpha_n} \right) =& 0 \\
\end{array}
\right. . $$ We know, using Assumption~\ref{hypforte} and the proof of Proposition~\ref{propasympnormphirstar2}, that $q_{\alpha_n}\left( R^* U^{(1)} \right)=\Phi_{R^*}^{-1}(\alpha_n)$ is asymptotically proportional to $(1-\alpha_n)^{-\frac{\gamma}{\gamma N+1}}$, while $A^* \left( \frac{1}{1-\alpha_n} \right)$ is asymptotically proportional to $\left( 1-\alpha_n \right)^{-\frac{\rho}{\gamma N+1}}$ if $\rho> -2 \gamma$ and $\left( 1-\alpha_n \right)^{\frac{2 \gamma}{\gamma N+1}}$ otherwise. Finally, it is not difficult to check that $\left( C_{int}^{L_p} \right)$ leads to the nullity of the two limits, and therefore to the third term of the decomposition, hence the result. The proof is exactly the same for the second normality, replacing $\hat{q}_{p,\alpha_n}\left(Y| \vec{X}=\vec{x} \right)$ by $\hat{\hat{q}}_{p,\alpha_n}\left(Y| \vec{X}=\vec{x} \right)$, $\ln \left( 1-\alpha_n \right)$ by $\ln \left( \frac{k_n}{n(1-\alpha_n)} \right)$ and using Proposition~\ref{propasympnormphirstar2} instead of~\ref{propnormint2}. $\Box$

\subsection*{Proof of Proposition~\ref{propnormhg}}

We have the following decomposition~:
\begin{multline}\nonumber
\frac{\sqrt{k_n}}{\ln \left( 1-\alpha_n \right)} \left( \frac{\hat{H}_{p,\alpha_n}\left(Y|\vec{X}=\vec{x} \right)}{H_{p,\alpha_n}\left(Y|\vec{X}=\vec{x} \right)}-1 \right) \underset{n \rightarrow +\infty}{\sim} \frac{\sqrt{k_n}}{\ln \left( 1-\alpha_n \right)} \left( \frac{\left( W_{\left[n \tilde{v}_n +1 \right]} \right)^{1/\hat{\eta}_{k_n}} f_H \left( \left( \hat{\gamma}_{k_n}^{-1}+N \right)^{-1}, p \right)}{H_{p,\alpha_n}\left(R^* U^{(1)} \right)}-1 \right) = \\ \frac{\sqrt{k_n}}{\ln \left( 1-\alpha_n \right)} \left( \frac{f_H \left( \left( \hat{\gamma}_{k_n}^{-1}+N \right)^{-1}, p \right)}{f_H \left( \left( \gamma^{-1}+N \right)^{-1}, p \right)}-1  \right) \frac{\left(W_{\left[n \tilde{v}_n +1 \right]} \right)^{1/\hat{\eta}_{k_n}}}{\Phi_{R^*}^{-1}(\alpha_n)} \frac{f_H \left( \left( \gamma^{-1}+N \right)^{-1}, p \right) \Phi_{R^*}^{-1}(\alpha_n)}{H_{p,\alpha_n}\left(R^* U^{(1)} \right)} + \\ \frac{\sqrt{k_n}}{\ln \left( 1-\alpha_n \right)} \left( \frac{\left(W_{\left[n \tilde{v}_n +1 \right]} \right)^{1/\hat{\eta}_{k_n}}}{\Phi_{R^*}^{-1}(\alpha_n)}-1 \right) \frac{f_H \left( \left( \gamma^{-1}+N \right)^{-1}, p \right) \Phi_{R^*}^{-1}(\alpha_n)}{H_{\alpha_n}\left(R^* U^{(1)} \right)}+ \\ \frac{\sqrt{k_n}}{\ln \left( 1-\alpha_n \right)} \left( \frac{f_H \left( \left( \gamma^{-1}+N \right)^{-1}, p \right) \Phi_{R^*}^{-1}(\alpha_n)}{H_{p,\alpha_n}\left(R^* U^{(1)} \right)}-1 \right).
\end{multline}
We know that $f_H \left( \left( \hat{\gamma}_{k_n}^{-1}+N \right)^{-1}, p \right)$, as a function of $\hat{\gamma}_{k_n}$, is asymptotically normal with rate $\sqrt{k_n}$ (see Equation~\eqref{eqnormgammahat}). Then, the first term in the sum clearly tends to 0 as $n \rightarrow +\infty$. Using Proposition~\ref{propnormint2}, the second term tends to the normal distribution given in~\eqref{eqnormint2}. Finally, we have to check that the third term tends to 0. For that purpose, we use the result of Theorem 4.5 in~\cite{maohu}, which ensures that there exists $\lambda \in \mathbb{R}$ such that~: $$ \frac{H_{p,\alpha_n}\left(R^* U^{(1)} \right)}{f_H \left( \left( \gamma^{-1}+N \right)^{-1}, p \right) \Phi_{R^*}^{-1}(\alpha_n)}=1+\lambda A^* \left( \frac{1}{1-\alpha_n} \right)(1+o(1)), $$ where $A^*$ is the auxiliary function of $\Phi_{R^*}^{-1} \left( 1-\frac{1}{t} \right)$. In the proof of Proposition~\ref{propnormint2}, we have seen that $A^*(t)$ was proportional either to $t^{-\frac{2 \gamma}{\gamma N+1}}$ if $\rho \leq -2 \gamma$ or $t^{\frac{\rho}{\gamma N+1}}$ otherwise. Then condition $\left( C_{int}^{HG} \right)$ ensures $$\underset{n \rightarrow +\infty}{\lim} \text{ } \frac{\sqrt{k_n}}{\ln \left( 1-\alpha_n \right)} \ln \left(  \frac{H_{\alpha_n}\left(R^* U^{(1)} \right)}{f_H \left( \left( \gamma^{-1}+N \right)^{-1}, p \right) \Phi_{R^*}^{-1}(\alpha_n)} \right) =0.$$ Hence the third term in the sum tends to 0, and the first result of~\eqref{eqnormhg} is proved. The proof is exactly the same for the second one, with rate $\frac{\sqrt{k_n}}{\ln \left( \frac{k_n}{n(1-\alpha_n)} \right)}$ instead of $\frac{\sqrt{k_n}}{\ln \left( 1-\alpha_n \right)}$. Then condition $\left( C_{high}^{HG} \right)$ gives the expected result.  $\Box$

\section*{Acknowledgements}

The author would like to thank the Editor-in-Chief, the Associate Editor and the referees, who did an extremely detailed and relevant report that widely helped to improve the paper. This work was partially supported by the MultiRisk LEFE-INSU Program, and by the LABEX MILYON (ANR-10-LABX-0070) of Universit\'e de Lyon, within the program ``Investissements d'Avenir'' (ANR-11-IDEX-0007) operated by the French National Research Agency (ANR).

\newpage

\bibliographystyle{apalike}
\bibliography{biblio2}

\end{document}